\documentclass[11pt,reqno]{amsart}

\usepackage{euscript}
\usepackage{pb-diagram}
\usepackage{lamsarrow}
\usepackage{pb-lams}
\usepackage{epsfig}
\usepackage{amssymb}

\theoremstyle{plain}

\theoremstyle{definition}

\newcommand{\no}{\noindent }
\newcommand{\bc}{{\mathbb C}}
\newcommand{\bp}{{\mathbb P}}
\newcommand{\bz}{{\mathbb Z}}

\newcommand{\bq}{{\mathbb Q}}

\begin{document}

\title{Links and analytic invariants  of superisolated singularities}
\author{I. Luengo-Velasco, A. Melle-Hern\'andez, and A. N\'emethi}
\address{Facultad de Matem\'aticas\\ Universidad Complutense\\ Plaza de Ciencias\\ E-28040, Madrid, Spain}
\address{Department of Mathematics\\Ohio State University\\Columbus, 
OH 43210,USA, \ and R\'enyi Institute of Math., Budapest, Hungary}
\email{iluengo@mat.ucm.es}
\email{amelle@mat.ucm.es}
\email{nemethi@math.ohio-state.edu, and nemethi@renyi.hu}
%\urladdr{http://www.math.ohio-state.edu/ \textasciitilde nemethi/}
\thanks{The first two authors are partially
supported by BFM2001-1488-C02-01. The third author is partially supported by NSF grant DMS-0304759.}



\begin{abstract} 
Using superisolated singularities we  present examples and counterexamples
to some of the most important conjectures regarding invariants 
of normal surface singularities. More precisely, we show that the 
``Seiberg-Witten invariant conjecture''(of Nicolaescu and the third author), 
the ``Universal abelian cover 
conjecture'' (of Neumann and Wahl)  and the ``Geometric genus conjecture'' 
fail (at least at that 
generality in which they were formulated). 
Moreover, we also show that for Gorenstein 
singularities (even with integral  homology sphere links) besides the 
geometric genus, the embedded dimension and the 
 multiplicity (in particular, the 
Hilbert-Samuel function)  also fail to be topological; and  
in general, the Artin cycle
does not coincide with the maximal (ideal) cycle.
\end{abstract}

\maketitle
\pagestyle{myheadings}
\markboth{{\normalsize I. Luengo-Velasco, A. Melle-Hern\'andez, A. 
N\'emethi}}{
{\normalsize Links and analytic invariants of SIS}}

{\small

\section{Introduction}

\subsection{}\label{1.1} In 
the last years we witness an intense effort to understand 
the following question: what kind of analytic invariants of an analytic
complex normal surface singularity can be determined from the topology 
(i.e. from the link) of the singularity ? Is the link indeed sufficiently
powerful to contain valuable information which would help to recover
analytic invariants (like multiplicity, Hilbert-Samuel function,
geometric genus), or equations (modulo equisingular deformations) ?
See, e.g. \cite{Artin62,Artin66,CoS,CoS3,FS,Laufer72,Laufer73,Laufer77,
Five,[36],[37],Ninv,NOSZ,Line,INV,
[51],[52],[55],Neu,NW,NWnew,NWnew2,NWuj,Ok,To1,
Yau4,Yau5,Yau1,Zariskiconj}. 
In fact, in order to give a chance to these type of questions, one has 
to assume two types of restrictions
(see e.g. \cite{Laufer77} and \cite{INV} for more details and examples): 
a topological one -- e.g. that the link is a rational homology sphere --
and an analytic one -- e.g. that the singularity is $\bq$--Gorenstein. 
Therefore, in the sequel we will assume that the link is a rational
homology sphere. 

As a result of the above mentioned efforts, in the last years a large number
of positive results and conjectures have appeared. Some of the conjectures 
were verified for large nontrivial families of singularities, a fact which 
created an increasing optimism. Nevertheless, some signs started to give
the signal that there are special families of singularities which might
create some obstructions, and whose understanding would be crucial for further
progress.
One of these families is the class of superisolated singularities.  

The goal of this note is to present examples and counterexamples
to some of the most important conjectures in this area regarding invariants 
of normal surface singularities (using 
superisolated singularities). More precisely, we show that the 
``Seiberg-Witten invariant conjecture''(of Nicolaescu and the third author), 
the ``Universal abelian cover 
conjecture'' (of Neumann and Wahl)  and the ``Geometric genus conjecture'' 
fail (at least at that 
generality in which they were formulated); see section 3 for a short review of
these conjectures. 
Moreover, these examples also show that for Gorenstein 
singularities (even with integral  homology sphere links) besides the 
geometric genus, the embedded dimension and the 
 multiplicity (in particular, the 
Hilbert-Samuel function)  also fail to be topological. 
One of the examples also shows that,  in general, the Artin cycle
does not coincide with the maximal (ideal) cycle (even for complete 
intersections with integral homology sphere links). 
The main message is 
that all the present conjectures and knowledge  must be reconsidered, 
rethought, reorganized in order to find the right and correct 
connections, directions and questions which would guide the next steps. 

Surprisingly, some of our examples are not very complicated (compared with 
the list of -- rather different and sometimes  rather complex --
positive examples which verify the 
corresponding conjectures). E.g., they are hypersurface singularities
(or their universal abelian covers). 
Nevertheless, they have some other rather specific properties
which allow some room for anomalies.

We wish to emphasize that the failure of the conjectures (at the generality
how they were formulated) puts in a different new light
all those families for which the conjectures were verified: 
their role and importance 
become much stronger and dominant. Moreover, this is a clear invitation for
clarification of some other new families of singularities, out of
which the superisolated singularities have the first priority. 

In section 2 we will set our notations and we will present some results
about the invariants of superisolated hypersurface singularities.
In section 3 we give the list of conjectures and problems for which 
we will provide examples-counterexamples in the following  sections. 
In section 4 our strategy is the following. We start with the classification
of the hypersurface superisolated singularities. By computing invariants
one gets directly counterexamples for the Seiberg-Witten invariant conjecture
(cf. 4.1) and the Universal abelian cover conjecture (cf. 4.3-4.4).
More complicated, but more striking examples are found by considering
the universal abelian cover of singularities (cf. 4.5-4.6).

\vspace{1mm}

The authors thank E. Artal-Bartolo, I. Dolgachev, J. Koll\'ar, J. Stevens
and J. Wahl for valuable discussions. 

\section{Hypersurface superisolated singularities}

\subsection{}\label{2.1} Hypersurface superisolated singularities
achieved historically the reputation of being an interesting class 
of  singularities.
This class ``contains'' in a canonical way the theory of complex projective 
plane curves, which gives a series of nice examples and counterexamples.
They were introduced in \cite{Ignacio} by the first author in order to show
that the $\mu$-constant stratum in the semiuniversal deformation space
of an isolated hypersurface singularity, in general, is not smooth.
Later Artal-Bartolo in \cite{Artal} used them to provide a counterexample
for S. S.-T. Yau's conjecture (showing that, in general, the link of
an isolated hypersurface surface 
singularity and its characteristic polynomial not determine the embedded topological type of the singular germ). 
On the other hand, A.~Durfee's conjecture and the 
monodromy conjecture of J.~Denef and F.~Loeser has been proved for them, 
see \cite{Melle} and \cite{aclm}.

\subsection{Definitions-Notations.}\label{2.2} 
A hypersurface singularity $f:(\bc^3,0)\to (\bc,0)$, $f=f_d+f_{d+1}+\cdots$
(where $f_j$ is homogeneous of degree $j$) is superisolated if the 
projective plane curve $C:=\{f_d=0\}\subset \bp^2$ is reduced with 
isolated singularities $\{p_i\}_i$, and these points are not situated 
on the projective curve $\{f_{d+1}=0\}$. In this case the embedded 
topological type (and the equisingular type) of $f$ does not depend on the 
choice of $f_{j}$'s (for $j>d$,
as long as $f_{d+1}$ satisfies the above requirement),
e.g.  one can take $f_j=0$ for any $j>d+1$ and $f_{d+1}=l^{d+1}$
where $l$ is a linear form not vanishing at the points
$\{p_i\}_i$. We will denote by $\mu_i$
(respectively by $\Delta_i$, with the sign choice $\Delta_i(1)=1$) 
the Milnor number (respectively, the characteristic polynomial)
 of the local plane curve singularities
$(C,p_i)\subset (\bp^2,p_i)$. For simplicity, in this note we will assume that
$C$ is irreducible. (The interested reader can adopt the next discussion 
easily to the general situation.)

Let $M$ be the link of $\{f=0\}$ (with its natural orientation),
$H:=H_1(M,\bz)$; $\mu$ 
and  $p_g$ be the Milnor number and the geometric genus of $f$. 

\subsection{Invariants.}\label{2.3} \

\vspace{1mm}

\no $\bullet$ \cite{Ignacio} The {\em  minimal 
resolution } of $\{f=0\}$ has only one irreducible exceptional divisor
which is isomorphic to $C$ and has self intersection $-d$.
In particular, the link $M$ of $f$ is a rational homology sphere if and 
only if $C$ is rational and all the plane curve singularities $(C,p_i)\subset 
(\bp^2,p_i)$ are (locally) irreducible (i.e., $C$ is a 
rational cuspidal plane curve).
In particular, $\sum_i\mu_i = (d-1)(d-2)$.

If $\Gamma_i$ is the minimal embedded resolution  graph of 
$(C,p_i)\subset  (\bp^2,p_i)$ (with a unique $-1$ vertex $v_i$
 which supports the 
strict transform of $(C,p_i)$), then the minimal good resolution graph
of $\{f=0\}$ can be constructed in the following way: 
consider a ``central vertex'' $v$
(which corresponds to the curve $C$), for each $i$ connect $v$ with $v_i$
by an edge, keep all the decorations of $\Gamma_i$, and add a new decoration 
$ e_v$ (self intersection) to $v $ as follows. In the graphs $\Gamma_i$ 
insert the set of multiplicities of the (reduced) plane curve singularity
(i.e. the strict transform of $(C,p_i)$ goes with multiplicity one).
Let $a_i$ be the multiplicity of the unique $-1$ curve of $\Gamma_i$.
Then $e_v=-d-\sum_ia_i$.

\vspace{1mm}

\no $\bullet$ Fix any resolution graph $\Gamma$ of $\{f=0\}$.
Let $K$ be the {\em canonical cycle} associated with $\Gamma$, and $s$ 
the number of vertices. Then $K^2+s$ does not depend on the choice of
$\Gamma$, it is a topological invariant of $M$.
In our case, it is easy to compute it at the level of the minimal
resolution:
$$K^2+s=-d(d-2)^2+1.$$

%\vspace{1mm}

\no $\bullet$ By (3.6.4) of \cite{Artal}, the {\em Milnor number }
$\mu$ of $f$ is the sum of the Milnor number of the singularity
$x^d+y^d+z^d$ and $\sum_i\mu_i$. Using the above mentioned identity
$\sum_i\mu_i=(d-1)(d-2)$, we get:
$$\mu=(d-1)^3+(d-1)(d-2).$$
Similarly, the {\em characteristic polynomial} $\Delta_f$ of $f$ is
$$\Delta_f(t)=\frac{t^d-1}{t-1}\cdot \prod_i\, \Delta_i(t^{d+1}).$$
Since $\Delta_i(1)=1$, this implies that $|H|=\Delta_f(1)=d$.
In fact, one can verify easily that $H=\bz_d$, and a possible generator
of $H$ is an elementary loop in a transversal slice to $C$. 

\vspace{1mm}

\no $\bullet$  Since $12p_g=\mu-(K^2+s)$ by \cite{Laufer77b},
one obtains:
$$p_g=d(d-1)(d-2)/6.$$

\vspace{1mm} 

\no $\bullet$ One of the conjectures relates the {\em 
Seiberg-Witten invariant} ${\bf sw}(M)$ of
$M$ (associated with the canonical $spin^c$ structure)
with the analytic (or smoothing) invariants of the singularity.
Here, by definition, ${\bf sw}(M)$ is
the sign-refined Reidemeister-Turaev torsion ${\mathcal T}(M)$
(associated with the canonical $spin^c$ structure)  \cite{Tu5} 
normalized by the Casson-Walker invariant $\lambda(M)$
 (using the convention of \cite{Lescop})
(cf. also with \cite{[51],[52],[55],NOSZ,INV}). Namely, we consider:
$${\bf sw}(M):=-\frac{\lambda(M)}{|H|}+{\mathcal T}(M).$$
Both invariants ${\mathcal T}(M)$ and $\lambda(M)$ can be determined from the 
graph (for details, see \cite{[51]} or \cite{INV}).
In fact, in our present case,
the formula of \cite{[51]} can be rewritten in the form:
$${\mathcal T}(M)=\frac{1}{d}\sum_{\xi^d=1\not=\xi}\ \frac{1}{(\xi-1)^2}
\cdot \prod_i\Delta_i(\xi).$$
Using similar method as in the proof of Theorem 4.5 of \cite{[55]}
(i.e. Fujita's splicing formula for the Casson-Walker invariant
\cite{fujita}, and Walker-Lescop surgery formula \cite{Lescop}, page 13)
one can proof the following identity.
Let $\bar{\Delta}(t)$ be the product 
$ \prod_i\, \Delta_i(t)$ symmetrized (i.e. its degree is 
$2\delta$ and $\bar{\Delta}(t)=t^{-\delta}\cdot \prod_i\, \Delta_i(t)$).
Then the Casson-Walker invariant of the link is 
$$\lambda(M)= (-1/2) \bar{\Delta}(t)''(1) + (d-1)(d-2)/24.$$
In fact, in this formula one can replace $\bar{\Delta}(t)''(1)$
by $\sum_i\bar{\Delta}_i(t)''(1)$. 

\section{The conjectures and questions}

\no Here we list the main conjectures and problems which have been
guiding our investigation.

\subsection{SWC. The Seiberg-Witten invariant conjecture.}\label{SWC} In 
\cite{[51]}
L. Nicolaescu and the third author formulated the following conjecture.

\vspace{2mm}

\no {\em (a) If the link of a normal surface singularity is a rational
homology sphere then 
$$p_g\leq {\bf sw}(M)-(K^2+s)/8.$$

\no (b) Additionally, if the singularity is ${\mathbb Q}$-Gorenstein, then
in (a) the equality holds.}

\vspace{2mm}

In the case of hypersurface singularities (more generally, for smoothings
of Gorenstein singularities) the identity (b) can be rewritten as 
$-8{\bf sw}(M)=\sigma$,
the signature of the Milnor fiber. If the singularity is  an isolated 
 complete intersection with an integral homology sphere link, then the 
conjecture transforms into the identity  $8\lambda(M)=\sigma$,
which was conjectured by Neumann and Wahl \cite{NW} for smoothings of 
complete intersections (this is called the ``Casson invariant conjecture'',
 $CIC$). 

The $SWC$-conjecture 
was verified e.g. for quotient singularities \cite{[51]},
for singularities with good 
$\bc^*$-actions \cite{[52]}, hypersurface suspension singularities
$g(u,v)+w^n$ with $g$ irreducible \cite{[55]}.
Even more, in \cite{NOSZ}, the third author replaced ${\bf sw}(M)$
by the corresponding Ozsv\'ath-Szab\'o invariant (which is 
defined via the Ozsv\'ath-Szab\'o Floer homology, and which 
conjecturally equals ${\bf sw}(M)$), and verified the inequality (a) for any
singularity with {\em almost rational} ({\em AR} in short)
resolution graph. 
(A graph is {\em AR}, if by replacing the decoration of 
one of the vertices one gets a rational graph. E.g., all the rational,
weakly elliptic, minimal good star-shaped graphs are {\em AR}.)

In fact, for rational singularities, even the equivariant
version of the $SWC$ was verified: \cite{Line} shows the identity
of the set of Seiberg-Witten invariants of the link (parametrized by
all the possible $spin^c$-structures) with the equivariant geometric 
genera of the universal abelian cover. 

On the other hand, in \cite{NOSZ} the author  exemplifies some types of graphs
which are not {\em AR}, and  whose understanding would be  necessary for 
further progress regarding the result of [loc.cit.].
These are exactly the type of graphs which are provided typically
by superisolated singularities. 

\subsection{UACC. The universal abelian cover conjecture.}\label{UACC}
The starting point of the next conjecture of Neumann and Wahl is 
 Neumann's paper \cite{Neu} which proves that the universal 
abelian cover of a singularity with a good $\bc^*$-action and with
$b_1(M)=0$ is a Brieskorn  complete intersection whose weights 
can be determined from the Seifert invariants of the link.
This, and other examples worked out by Neumann and Wahl (see e.g.
\cite{NWnew2} about quotient-cusps)  lead them to 
a rather complex program and package of conjectures \cite{NWnew}:

\vspace{2mm}

\no {\em Assume that $(X,0)$  is $\bq$-Gorenstein 
singularity with $b_1(M)=0$. 
Then there exists an equisingular and equivariant 
deformation of the universal  abelian cover of $(X,0)$ 
to an isolated  complete intersection singularity. Moreover, 
the equations of this complete intersection, 
together with the action of $H_1(M,\bz)$,
can be recovered from $M$ via the ``splice equations''}. 

\vspace{2mm}

The main point of the above conjecture, in its detailed
version, provides a clear recipe for the equations
of the complete intersection singularity (the ``splice equations'')
and the action of $H$ on these equations. This is done in terms of the 
combinatorics of the resolution graph of $(X,0)$.  
In order to be able to write
down the equations, the graph should satisfy some 
arithmetical properties: the so-called {\em semigroup conditions}
and {\em congruence conditions}. Their validity is part of the 
conjecture. The reader is invited to see all the details
in \cite{NWnew}. 

In order to eliminate any confusion, we mention that in this note
 {\em equisingular deformation} means the existence of a 
{\em simultaneous equitopological resolution} as discussed in 
\cite{WahlAnnals}. 

\subsection{GGC. The geometric genus conjecture.}\label{GGC} Both
\ref{SWC} and \ref{UACC} are closely related with the following 
more general conjecture, which was formulated as a very general 
guiding principle
(cf. with Question (3.2) in \cite{NW}, see also  Problem 9.2 in \cite{INV}). 

\vspace{2mm}

\no {\em In the case of a  $\bq$-Gorenstein singularity with $b_1(M)=0$,
the geometric genus $p_g$ is topological (i.e. can be recovered from
the oriented link).}

\vspace{2mm}

Here we mention the following positive result of Pinkham \cite{Pi1}:
If a singularity with $b_1(M)=0$ has a good $\bc^*$-action, then its geometric
genus can be computed explicitly from the resolution graph.
(Moreover, by \cite{Neu}, such a singularity is $\bq$-Gorenstein.)

The crucial testing case for the above $GGC$ is the case of the 
star-shaped resolution graphs: {\em is it true that if the minimal good
resolution graph of a $\bq$-Gorenstein singularity is star-shaped, then its 
geometric genus is the same as the number predicted by Pinkham's formula ?}

\subsection{Other analytic invariants.}\label{other} Similar 
questions were raised for several other discrete analytic invariants
as well:

\vspace{2mm} 

\no {\em For what family of $\bq$-Gorenstein singularities (with $b_1(M)=0$)
are the invariants like the multiplicity, embedded dimension,
Hilbert-Samuel function, maximal cycle (etc.) topological ?}

\vspace{2mm}

For different positive cases and comments, see  e.g. \cite{INV}.
(The fact that the embedded dimension can jump 
in a topological constant family  -- even in a 
positive-weight deformation of a weighted homogeneous singularity -- 
was known by experts.)

\vspace{2mm}

{\em
The examples of the next sections provide negative answers to all of the above
conjectures ($SWC$, $UACC$ and $GGC$)
and all the analytic invariants   listed in \ref{other}
(already in Gorenstein case).}

\section{Examples}
 
\subsection{Examples/counterexamples for the $SWC$-conjecture.}\label{4.1} \ 
Some of the next examples show that the $SWC$-conjecture, 
in general, is not true. 
For this, we consider superisolated singularities  $f=f_d+l^{d+1}$.
Below, any singular point $(C,p_i)$ will be identified by its multiplicity 
sequence. Since the number of occurrences of the multiplicity $1$
in the multiplicity sequence equals the last multiplicity greater than $1$,
we omit the multiplicity $1$: we denote such a sequence by
$ [m_0,\ldots,m_l]$ where $ m_0 \ge m_1 \ge \cdots \ge m_l > m_{l+1} = 1$
for a suitable $l \ge 0$. In fact, we 
will write $[\hat{m}_{0_{r_0}},\ldots,\hat{m}_{k_{r_k}}]$ 
for a multiplicity sequence which means that the multiplicity $\hat{m}_i$ occurs $r_i$ times for
$i=0,\ldots,k$. For example, $[4_2,2_3]$ means $[4,4,2,2,2,1,1]$.

If $C$ has only one singularity with sequence $[d-1]$, 
then $f$ has the same invariants as 
the weighted homogeneous singularity $zx^{d-1}+y^d +z^{d+1}$, hence
it satisfies the conjecture by \cite{[52]}. 
(Probably it is worth to mention that not all the rational cuspidal 
curves of degree $d$ with one cusp and multiplicity sequence $[d-1]$
are projectively equivalent. E.g., for $d=4$ there are two 
projectively non-equivalent curves: $\{x^4-x^3y+y^3z=0\}$  and
$\{x^4-y^3z=0\}$, cf. \cite{Namba}, page 135.) 

If $d=3$, then $C$ has a unique singularity of type $[2]$. 
If $d=4$, then there are four possibilities; the corresponding 
multiplicity sequences of the singular points $\{p_i\}_i$ of $C$ are 
$[3]$; $[2_3]$; $[2_2],[2]$ and $[2],[2],[2]$. 
By a verification, in all  these cases, the conjecture is again true.
(For the classification of the cuspidal rational curves with small
degree, see e.g. the book of Namba \cite{Namba}.)

If $d=5$, then $p_g=10$, $K^2+s=-44$. Let $N$ be the number of singular 
points of $C$. The next table shows for all the possible multiplicity
sequences 
the validity of the conjecture. When the conjecture fails, we put in 
parenthesis the value $-\lambda/|H|+{\mathcal T}-(K^2+s)/8$
(which can be compared with the value of $p_g$).

\vspace{2mm}

\begin{enumerate}
  \item $N=1$

  \vspace{0.0cm}
    \hspace{8ex}
    \begin{tabular}{c|l|c}
              & type of cusp       & $conj$ \\ \hline 
      $C_1$   & $[4]$           & $True$         \\ 
      $C_2$   & $[2_6]$      & $True$        \\
     
    \end{tabular}

   \vspace{0.0cm}
  \item $N=2$
 
   \hspace{8ex}
    \begin{tabular}{c|l|c} 
              & type of cusps & $conj$ \\ \hline 
      $C_3$   & $[3,2]$ , $[2_2]$  & $True$  \\ 
      $C_4$   & $[3]$ , $[2_3]$  &  $False$ (8)  \\ 
      $C_5$   & $[2_2]$ , $[2_4]$ & $False$ (8)\\
      
    \end{tabular}

    \vspace{0.0cm}

  \item $N=3$

    \hspace{7.2ex}
    \begin{tabular}{c|l|c}
            & type of cusps & $conj$ \\ \hline  
      $C_{6}$   & $[3]$ , $[2_2]$ , $[2]$ & $False$ (8)   \\ 
      $C_{7}$   & $[2_2]$ , $[2_2]$ , $[2_2]$ & $False$ (4)  \\ 
    \end{tabular}
  
    \vspace{0.0cm}

  \item $N=4$

    \hspace{7.2ex}
    \begin{tabular}{c|l|c}
            & type of cusps & $conj$ \\ \hline  
      $C_{8}$   & $[2_3]$ , $[2]$ ,$[2]$, $[2]$& $False$ (2)   \\ 
      
    \end{tabular}
    \end{enumerate}

\vspace{1mm}

\no If $d=6$ then $p_g=20$. 
The classification of multiplicity sequences of 
rational cuspidal plane curves of degree $6$ with $N$ singular points
is given by the following list, see e.g.  Fenske's paper \cite{Fenske}.

\vspace{1mm}

\begin{enumerate}
  \item $N=1$

  \vspace{0.0cm}
    \hspace{8ex}
    \begin{tabular}{c|l|c}
              & type of cusp       & $conj$ \\ \hline 
      $C_1$   & $[5]$           & $True$         \\ 
      $C_2$   & $[4,2_4]$      & $True$        \\
      $C_3$   & $[3_3,2]$      &  $True$        \\
    \end{tabular}
    
   \vspace{0.0cm}

  \item $N=2$
 
   \hspace{8ex}
    \begin{tabular}{c|l|c} 
              & type of cusps & $conj$ \\ \hline 
      $C_4$   & $[3_3]$  , $[2]$  & $True$  \\ 
      $C_5$   & $[3_2,2]$ , $[3]$  &  $True$  \\ 
      $C_6$   & $[3_2]$      , $[3,2]$ & $True$ \\
      $C_7$   & $[4,2_3]$ , $[2]$  & $True$  \\ 
      $C_8$   & $[4,2_2]$ , $[2_2]$  &  $True$ \\ 
      $C_9$   & $[4]$      , $[2_4]$  &  $False$ (18)  \\
    \end{tabular}

    \vspace{0.0cm}

  \item $N=3$

    \hspace{7.2ex}
    \begin{tabular}{c|l|c}
            & type of cusps & $conj$ \\ \hline  
      $C_{10}$   & $[4]$ , $[2_3]$ , $[2]$ & $True$   \\ 
      $C_{11}$   & $[4]$ , $[2_2]$ , $[2_2]$ & $True$   \\ 
    \end{tabular}
  
    \vspace{0.0cm}
  \end{enumerate}

\vspace{1mm}

\no For the convenience of the reader, we make the example $d=5$, $N=2$,
case $C_4$, more explicit.  In this case the minimal good resolution graph
has the form

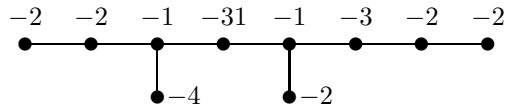
\begin{picture}(300,45)(0,0)
\put(100,25){\circle*{5}}
\put(125,25){\circle*{5}}
\put(150,25){\circle*{5}}
\put(175,25){\circle*{5}}
\put(200,25){\circle*{5}}
\put(225,25){\circle*{5}}
\put(250,25){\circle*{5}}
\put(275,25){\circle*{5}}
\put(150,5){\circle*{5}}
\put(200,5){\circle*{5}}
\put(100,25){\line(1,0){175}}
\put(150,25){\line(0,-1){20}}
\put(200,25){\line(0,-1){20}}
\put(100,35){\makebox(0,0){$-2$}}
\put(125,35){\makebox(0,0){$-2$}}
\put(150,35){\makebox(0,0){$-1$}}
\put(175,35){\makebox(0,0){$-31$}}
\put(200,35){\makebox(0,0){$-1$}}
\put(225,35){\makebox(0,0){$-3$}}
\put(250,35){\makebox(0,0){$-2$}}
\put(275,35){\makebox(0,0){$-2$}}
\put(160,5){\makebox(0,0){$-4$}}
\put(210,5){\makebox(0,0){$-2$}}
\end{picture}

\no By \cite{[51]}, or by the above formulae,
$-\lambda(M)=21/2$ and ${\mathcal T}(M)=2/5$,
hence ${\bf sw}(M)-(K^2+s)/8=8$. 

\vspace{2mm}

Notice that above, in all the cases when part (b) of the 
$SWC$-conjecture fails
(i.e. $p_g\not={\bf sw}(M)-(K^2+s)/8$), part (a) of \ref{SWC}  fails as
 well: the topological candidate becomes strict {\em smaller} than $p_g$.

%Notice that none of these examples contradicts the most general 
%guiding conjecture, which claims  that $p_g$ can be recovered from 
%the link (provided that the link is a rational homology sphere and 
%the singularity is ${\mathbb Q}$-Gorenstein). In order to produce a 
%counterexample for this last statement one has to produce a pair 
%of singularities with the same graph but different $p_g$.
%In our present situation we are in some sense in the opposite
%situation: the Seiberg-Witten invariant is too sensible to the 
%variation of the graph (it probably keeps too much
%information from the graph). 

\subsubsection{}\label{4.1b} The authors analyzed even higher degree
 curves $C$ present in the literature, 
but were not able to find any counterexample with $N=1$. 
Although this very paper shows how cautious one should be with formulating
conjectures, still,  
 {\em we predict that for $N=1$ the $SWC$ is actually true}.
This conjecture  is also supported by its verification for a series
of non-trivial families, e.g. for irreducible curves
$C$ of {\em Abhyankar-Moh-Suzuki type}. They are characterized by the 
existence of  a line $L\subset \bp^2$ such that $C\setminus L
$ is isomorphic to $\bc$ (or, $C\cap L$ is the unique singular point
of $C$). (Notice that 
not any curve with $N=1$ satisfies this property, e.g. the Yoshihara quintic 
-- $C_2$ with $[2_6]$ in our table in \ref{4.1} -- does not.)
We also verified the above conjecture for all the cases when 
the singular point has exactly one characteristic pair. 
Since the techniques involved in these  verifications are rather different
from the spirit of the present note, they will be presented in another 
article \cite{4authors}. 

In fact, since any hypersurface superisolated singularity with $N=1$
is {\em AR} (in the sense of \cite{NOSZ}), the inequality 
\ref{SWC}(a) is valid for them by \cite{NOSZ}, 9.5(a).

\subsection{Remark.} \ 
Analyzing the above examples \ref{4.1}, 
one can ask: why the class of superisolated
singularities is so special?  In the spirit of \cite{NOSZ} (i.e.
thinking about non-{\em AR} graphs) we can notice that if we want to 
transform the above superisolated graphs into rational graphs by replacing
the original self intersection numbers by more negative ones, 
then we have to do this for many vertices (at least for $N$ vertices 
 of type $v_i$). Is the presence of these ``bad'' vertices the reason
for the above anomalies? The answer probably is that not just this:
one can produce easily 
 suspension singularities (which verify the conjecture by 
\cite{[55]}) with more than one ``bad'' vertex. For example, let $g(x,y)$ be 
the  irreducible plane curve singularity with Newton pairs $(p_1,q_1)=
(5,6)$ and $(p_2,q_2)=(2,5)$. Then the resolution graph of the suspension
hypersurface singularity $f(x,y,z)=g(x,y)+z^5$ is the following
(for the corresponding algorithm, see \cite{Five}, Appendix):

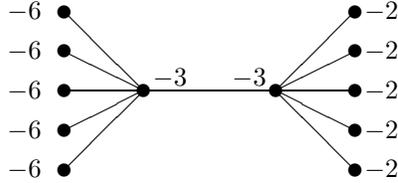
\begin{picture}(300,80)(0,-10)
\put(100,30){\circle*{5}}
\put(70,30){\circle*{5}}
\put(150,30){\circle*{5}}
\put(70,45){\circle*{5}}
\put(70,60){\circle*{5}}
\put(70,15){\circle*{5}}
\put(70,0){\circle*{5}}
\put(180,30){\circle*{5}}
\put(180,45){\circle*{5}}
\put(180,60){\circle*{5}}
\put(180,15){\circle*{5}}
\put(180,0){\circle*{5}}
\put(100,30){\line(1,0){50}}
\put(100,30){\line(-2,1){30}}
\put(100,30){\line(-1,1){30}}
\put(100,30){\line(-1,0){30}}
\put(100,30){\line(-2,-1){30}}
\put(100,30){\line(-1,-1){30}}
\put(150,30){\line(2,1){30}}
\put(150,30){\line(1,1){30}}
\put(150,30){\line(1,0){30}}
\put(150,30){\line(2,-1){30}}
\put(150,30){\line(1,-1){30}}
\put(55,0){\makebox(0,0){$-6$}}
\put(55,15){\makebox(0,0){$-6$}}
\put(55,30){\makebox(0,0){$-6$}}
\put(55,45){\makebox(0,0){$-6$}}
\put(55,60){\makebox(0,0){$-6$}}
\put(190,0){\makebox(0,0){$-2$}}
\put(190,15){\makebox(0,0){$-2$}}
\put(190,30){\makebox(0,0){$-2$}}
\put(190,45){\makebox(0,0){$-2$}}
\put(190,60){\makebox(0,0){$-2$}}
\put(110,35){\makebox(0,0){$-3$}}
\put(140,35){\makebox(0,0){$-3$}}
\end{picture}

\noindent 
In this case the graph has two ``bad''  vertices, 
$H=\bz_6^4\oplus \bz_2^4$, $K^2+s=-244$, $\mu=416$,
$p_g=55$, $-\lambda/|H|=61/18$, ${\mathcal T}=190/9$, and 
${\bf sw}=49/2=-\sigma/8$ (for different formulas and details
regarding suspension singularities, see e.g. \cite{[55]}). 

\subsection{Counterexamples for $UACC$.}\label{4.2} \ 
Working with superisolated singularities one sees easily that
already the construction of the ``splicing equations'' 
is obstructed: in general, the {\em semigroup condition is not 
satisfied}. More precisely, consider the splice diagram
associated with the resolution graph of a hypersurface superisolated 
singularity. Then, if $N\geq 3$, that decoration
of any edge of type $[v,v_i]$ which is closer to $v$ is 1. 
This should be situated in the semigroup generated by the 
decorations of the leaves (which are all strict greater than 1), a fact
which is not true. 

{\em This means that the algorithm \cite{NWnew}
which provides the equations of the 
complete intersection singularity predicted by the $UACC$
is not working, since to write the splice complete intersection
equations one needs the semigroup condition satisfied.
In other words, that conjectured complete intersection singularity whose
deformation should contain the universal abelian cover,
in general, does not exist.}

(But, of course, this does not imply that the universal abelian  cover cannot 
be a complete intersection; it might be, but not of splice type.)

The simplest counterexample appears when $d=4$ and $C$ 
is the Steiner quartic (the unique quartic 
in the plane with three $[2]$-cusps).

\subsection{A suspension type counterexample.}\label{4.3} \ 
In fact, the phenomenon \ref{4.2} is not really specific to superisolated 
singularities. One can construct hypersurface suspension singularities
with the same property. 
E.g. if one takes the hypersurface singularity
$\{z^2=(y+x^2)(y^3+x^{11})\}$, then its link is a rational homology sphere
(with first homology $\bz_4$), but its minimal plumbing  graph does not 
satisfy the semigroup  conditions (since the $E_8$-subgraph has
determinant 1). The resolution graph is 

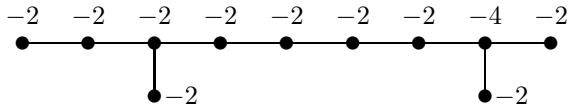
\begin{picture}(300,45)(0,0)
\put(50,35){\makebox(0,0){$-2$}}
\put(75,35){\makebox(0,0){$-2$}}
\put(100,35){\makebox(0,0){$-2$}}
\put(125,35){\makebox(0,0){$-2$}}
\put(150,35){\makebox(0,0){$-2$}}
\put(175,35){\makebox(0,0){$-2$}}
\put(200,35){\makebox(0,0){$-2$}}
\put(225,35){\makebox(0,0){$-4$}}
\put(250,35){\makebox(0,0){$-2$}}
\put(110,5){\makebox(0,0){$-2$}}
\put(235,5){\makebox(0,0){$-2$}}
\put(50,25){\circle*{5}}
\put(75,25){\circle*{5}}
\put(100,25){\circle*{5}}
\put(125,25){\circle*{5}}
\put(150,25){\circle*{5}}
\put(175,25){\circle*{5}}
\put(200,25){\circle*{5}}
\put(225,25){\circle*{5}}
\put(250,25){\circle*{5}}
\put(100,5){\circle*{5}}
\put(225,5){\circle*{5}}
\put(50,25){\line(1,0){200}}
\put(100,5){\line(0,1){20}}
\put(225,5){\line(0,1){20}}
\end{picture}

\vspace{3mm}

More sophisticated counterexamples are provided by considering universal 
abelian covers. 

\vspace{3mm}

\subsection{Counterexample:
The case $d=4$ with multiplicity sequence $[2_3]$, and 
its universal abelian cover.}\label{4.4} \
In this section, we make more explicit the invariants of the superisolated 
singularity with $d=4$ when $C$ has only one singular point of type $[2_3]$.
In this case $C$ is projectively equivalent to the projective curve
$(zy-x^2)^2=xy^3$, with parametrization $[t:s]\mapsto 
[t^2s^2: t^4: s^4+t^3s]$ (see \cite{Namba}, page 146). 
Hence a possible choice for $f$ is
$$f=(zy-x^2)^2-xy^3+z^5.$$
As we already mentioned $p_g(X,0)=4$. 
The resolution graph  $\Gamma$ is 

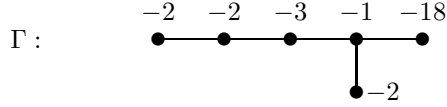
\begin{picture}(300,45)(0,0)
\put(100,25){\circle*{5}}
\put(125,25){\circle*{5}}
\put(150,25){\circle*{5}}
\put(175,25){\circle*{5}}
\put(200,25){\circle*{5}}
\put(175,5){\circle*{5}}
\put(100,25){\line(1,0){100}}
\put(175,25){\line(0,-1){20}}
\put(50,25){\makebox(0,0){$\Gamma:$}}
\put(100,35){\makebox(0,0){$-2$}}
\put(125,35){\makebox(0,0){$-2$}}
\put(150,35){\makebox(0,0){$-3$}}
\put(175,35){\makebox(0,0){$-1$}}
\put(200,35){\makebox(0,0){$-18$}}
\put(185,5){\makebox(0,0){$-2$}}
\end{picture}

\noindent 
Since the graph is a star-shaped, the same resolution graph
can be realized by a weighted homogeneous singularity $(X_w,0)$ as well.
In the present case this is an isolated
 complete intersection in $(\bc^4,0)$ with two equations:
$$ (X_w,0)= \left\{ \begin{array}{l}
yz=x^2\\
z^5+t^2-xy^3=0.\end{array}\right.$$
The corresponding weights of the coordinates $(x,y,z,t)$ are: 
 %8/35, 9/35, 1/5 and 1/2.
(16,18,14,35). 
 By Pinkham's formula \cite{Pi1} one gets that
$p_g(X_w,0)=4$ as well. 

Now, one can use the result of Neumann and Wahl (cf. \cite{NW} (3.3))
which guarantees that a  Gorenstein singularity (with the same link as 
a weighted homogeneous singularity $(X_w,0)$) is an 
equisingular deformation of $(X_w,0)$ if and only
if its $p_g$ equals the number predicted by Pinkham's formula.
(In \cite{NW} (3.3) the result is stated for integral homology spheres links,
but the proof works without modification for rational homology sphere 
links as well.)

In particular, the superisolated singularity $(X,0)$ is an
equisingular deformation of $(X_w,0)$. In this particular case
this deformation can be written easily: the pair of equations 
$yz-x^2=\lambda t, $ and $z^5+t^2-xy^3=0$ 
-- for  the parameter $\lambda\not=0$ -- is equivalent to $(X,0)$. 
(Here, if one wishes to emphasize the compatibility of the weights
with the deformation, one should notice that 
$\lambda$ has weight $-3$).

\vspace{2mm}

Notice also that the two singularities $(X,0)$ and $(X_w,0)$
{\em have the same multiplicity} (which equals 4), but clearly have 
{\em different
embedded dimensions -- hence different Hilbert-Samuel functions}. 

\vspace{2mm}

Next, we wish to analyze the corresponding universal abelian covers.

The universal abelian cover $(X_w^{ab},0)$ of $(X_w,0)$ is easy
(cf. also with \cite{Neu}). 
It is a hypersurface Brieskorn singularity $\{u^7=v^{18}+w^2\}$.
(Notice that this equation is exactly the ``splice equation'' predicted 
by the Neumann-Wahl construction, cf. \ref{UACC}). 
The action of (a generator of) $H=\bz_4=\{\zeta\in \bc\,:\, \zeta^4=1\}$ 
on the coordinates $(u,v,w)$ is
$u\mapsto -u$, $v\mapsto iv$ and $w\mapsto -iw$. 
Taking the invarints $x:=uv^2$, $y:= u^2$, $z:=v^4$ and $t:=vw$, we get 
that $\{u^7=v^{18}+w^2\}/H$ has exactly those equations what we provided 
for $(X_w,0)$ above. 

Notice that the resolution graphs of both universal abelian covers are
the same (which is exactly the plumbing diagram of the universal abelian cover
of the common link $M$). In this case it is:

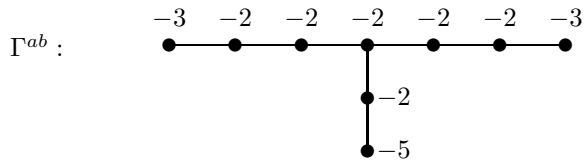
\begin{picture}(300,65)(0,-20)
\put(50,25){\makebox(0,0){$\Gamma^{ab}:$}}

\put(100,25){\circle*{5}}
\put(125,25){\circle*{5}}
\put(150,25){\circle*{5}}
\put(175,25){\circle*{5}}
\put(200,25){\circle*{5}}
\put(225,25){\circle*{5}}
\put(250,25){\circle*{5}}

\put(175,5){\circle*{5}}
\put(175,-15){\circle*{5}}

\put(100,25){\line(1,0){150}}
\put(175,25){\line(0,-1){40}}

\put(100,35){\makebox(0,0){$-3$}}
\put(125,35){\makebox(0,0){$-2$}}
\put(150,35){\makebox(0,0){$-2$}}
\put(175,35){\makebox(0,0){$-2$}}
\put(200,35){\makebox(0,0){$-2$}}
\put(225,35){\makebox(0,0){$-2$}}
\put(250,35){\makebox(0,0){$-3$}}
\put(185,5){\makebox(0,0){$-2$}}
\put(185,-15){\makebox(0,0){$-5$}}
\end{picture}

This graph has $K^2+s=-18$. Since the Brieskorn singularity $u^7=v^{18}+w^2$
has Milnor number $6\cdot 17=102$, we get by Laufer's formula that
its geometric genus is $p_g(-u^7+v^{18}+w^2)=10$.

Next, we analyze the universal abelian cover of the superisolated
singularity $(X,0)$ and estimate  its geometric genus.

In our original study of examples of \ref{4.4} and \ref{4.5},
the authors  had the faulty impression that the equisingular deformation 
existing  at the level of $(X,0)$ lifts to an equisingular and
equivariant deformation at the level of the universal abelian cover.
But when we showed J. Wahl the example \ref{4.5}, he recognized that this
could indeed not occur, and outlined a proof: any equivariant positive weight
deformation of the universal abelian cover of $(X_w,0)$ gives a family
of quotients of constant embedding dimension 6. 

In the sequel we present an alternate proof of the non-existence of
such deformation (using Fact B below, which is interesting by its own,
and hopefully can be applied in different similar situations as well).

\vspace{2mm}

\no {\bf Fact A.} {\em The universal abelian cover $(X^{ab},0)$ of
$(X,0)$ is not in the $\mu$-constant deformation space of
$(X_w^{ab},0)=\{w^2+v^{18}-u^7\}$. In particular, there is no equisingular
deformation from $(X^{ab},0)$ to $(X_w^{ab},0)$.}

\vspace{2mm}

In fact, what we will prove is the following: 

\vspace{2mm}

\no {\bf Fact B.} {\em Assume that $\bz_4$ acts
freely in codimension 1 on a hypersurface germ which in some coordinates
has the form $w^2+(\deg \geq 5)$. Then if the quotient is a 
hypersurface with multiplicity greater than $2$,
then the tangent cone  of the quotient is reducible. }

\begin{proof}
\no (1) Notice that any germ in the semiuniversal deformation of $(X_w^{ab},0)$
(modulo a coordinate change) can be written in the form $w^2+g(u,v)$.
Assume that $(X^{ab},0)$ is given by $(\{f^{ab}=0\},0)\subset (\bc^3,0)$
and $f$ is in a $\mu$-constant deformation of $(X_w^{ab},0)$.
%Notice that any germ in the semiuniversal deformation of $(X_w^{ab},0)$
%(modulo a coordinate change) can be written in the form $w^2+g(u,v)$.
Hence $f^{ab}$ itself, in some coordinates, has this form such that 
the plane curve singularities $u^7-v^{18}$ and $g(u,v)$
have the same embedded topological types. 
Therefore, all the monomials of $g$ have degree at least 7.

\no (2) We consider the action of $\bz_4 $ on $\{f^{ab}=0\}$. Since $\{f^{ab}=0\}$
is singular with tangent space $\bc^3$, the action induces an action
on this cotangent space $m/m^2$ and on the exact sequence $0\to m^2\to m
\to m/m^2\to 0$ (here $m$ is the maximal ideal of $\bc[u,v,w]/(f^{ab})$).
Since $\bz_4$ is finite, this sequence equivariantly splits, hence the 
singularity has an equivariant embedding  into its tangent space. In other 
words, by a change of local coordinates, we can assume that that the action 
extends to a linear action of $\bc^3$. Since the group is cyclic 
(with distinguished generator $\epsilon$), we can even assume that the
linear action on $\bc^3$ is diagonal. 

Since the space $\{f^{ab}=0\}$ is invariant to the action, $f^{ab}$ is an eigenfunction
of $\epsilon$ (coinvariant). Since $f^{ab}$ (in any coordinates) has the form
$l^2+(\deg\geq 3)$, where $l$ is a linear form, $l^2$ is also an eigenfunction
 of $\epsilon$ with the same eigenvalue $\lambda$
as $f^{ab}$. Since $l^2$ is a square, we get that $\lambda=\pm1$. 
Moreover, if $l$ involves more coordinates with nonzero coefficients,
then the action of $\epsilon $ on all of them should be the same, hence by
another linear change of  variables, and keeping the diagonal form of 
$\epsilon$, we can transform $l$ into one of the coordinates. 
We will denote the  coordinates constructed in this way  by $w,u,v$.

The action of $\epsilon$ has the form $diag(i^{a_1},i^{a_2},i^{a_3})$.
Since the action on $\{f^{ab}=0\}$ is free in codimension 1, one gets $\#
\{j:\, a_j \  \mbox{even}\}\leq 1$. 

\no (3) Consider the projection $p:\bc^3 \to \bc^3/\bz_4$.
If $g$ is a function vanishing along $\{f^{ab}=0\}/\bz_4$, then 
$g\circ p$ is an invariant function of form $f^{ab}h$. 
In particular, in order to obtain all the equations of $\{f^{ab}=0\}/\bz_4$,
we have to multiply $f^{ab}$ with such coinvariants $h$ which make 
$f^{ab}h$ invariant, and  express $f^{ab}h$ in terms of principal invariants.
If $f^{ab}$ itself is invariant, it provides basically only one equation, namely
its expression in terms of the principal invariants.

\no (4) Assume that there is an $a_j$ (say $a_1$) multiple of 4. 
Then (modulo some symmetry) there are the following possibilities:

\no (4.1) $\epsilon=diag(1,i,i)$ (i.e. $\epsilon(w)=w$,
$\epsilon(u)=iu$, $\epsilon(v)=iv$). The 
principal invariants in $\bc [w,u,v]$ of the action are
$I=\{w,v^4,v^3u,v^2u^2,vu^3,u^4\}$, hence $embdim(\bc^3/\bz_4)=6$.
Recall that $f^{ab}=l^2+(\deg\geq 3)$, where $l$ is one of the coordinates.

\no (4.1.1) Assume that $f^{ab}=w^2+(\deg\geq 3)$. Since $f^{ab}$ in some coordinates
has the from $\bar{w}^2+(\deg\geq 5)$ (cf. part (1)),
$f^{ab}$ in variables $(w,u,v)$ can be written as
\begin{equation*}
f^{ab}=(w+h_2+h_3)^2+(\deg\geq 5)=w^2+2wh_2+h_2^2+2wh_3+(\deg\geq 5),
\tag{$*$}
\end{equation*}
where $\deg(h_j)=j$. In this case $f^{ab}$ and $w$ are invariants, hence the 
same is valid for $wh_2$ and $h_2^2+2wh_3$ as well. Hence $f^{ab}=w^2+aw^3+bw^4
+(\deg\geq 5)$; in particular $f^{ab}$ expressed in terms of the invariants $I$
has no linear term. Therefore, $embdim\{f^{ab}=0\}/\bz_4=6$.

\no (4.1.2) Assume that $f^{ab}=u^2+(\deg\geq 3)$ (the case $f^{ab}=v^2+\cdots$ is 
symmetric). Then the principal invariants  $u^4,vu^3,v^2u^2$ can be 
eliminated using the equation of $f^{ab}$. The remaining relevant principal
invariants are $x:=w,y:=v^4,z:=v^3u$. Hence $\{f^{ab}=0\}/\bz_4$ can be embedded
into $(\bc^3,0)$. Next we analyze its equation. Notice that in this case 
$\lambda=-1$. Therefore, if $m=w^\alpha u^\beta v^\gamma$ is a monomial of $f^{ab}$
with nonzero coefficient, then $\beta+\gamma=4t_m+2$ for some $t_m\geq 0$.
If we multiply this monomial by $v^{4k+6}$, we get the invariant 
$w^\alpha u^\beta v^{\gamma+4k+6}$. Notice that the inequality
$\gamma+4k+6\geq 3\beta$ is equivalent with $k\geq \beta-t_m-2$, hence if we 
take  $k_0:=\max_m(\beta-t_m-2)$, then 
$mv^{4k_0+6}=x^\alpha z^\beta y^{t_m-\beta +2+k_0}$. 
In particular, $u^2v^{4k_0+6}=z^2y^{k_0}$. 
In other words, if $f^{ab}=\sum_ma_mm$, then the wanted equation of the 
quotient in $(\bc^3,0)$ is $\sum_ma_m x^\alpha z^\beta y^{t_m-\beta +2+k_0}$. 

Finally notice that if the $(w,u,v)$-degree of $m$ is
$\alpha+\beta+\gamma>2$ (i.e. if $m$ is any monomial different from
$u^2$), then the  $(x,y,z)$-degree of $mv^{4k_0+6}$
 is $\alpha+\beta+t_m-\beta+2+k_0>2+k_0$.
This shows that $\{f^{ab}=0\}/\bz_4$ in $(\bc^3,0)$ (with coordinates 
$x,y,z$)
has the equation $z^2y^{k_0}+$ higher degree terms. 
This, by any coordinate change, is not equivalent with the superisolated 
hypersurface singularity $f_4+f_5$ (because its tangent cone is 
reducible).

\no (4.2) Assume that $\epsilon$ acts on $(w,u,v)$
diagonally via $diag(1,i,-i)$. The set of principal invariants are 
$I=\{w,v^4,uv,u^4\}$. 

\no (4.2.1) If $f^{ab}=w^2+\ldots$, then using the same notation as in ($*$),
$h_2$ is invariant, hence $f^{ab}$, expressed in terms of the principal
invariants, has no linear term. In particular, $embdim\{f^{ab}=0\}/\bz_4=4$.

\no (4.2.2) If $f^{ab}=u^2+(\deg\geq 3)$, we proceed as in (4.1.2).
The invariant $u^4$ can be eliminated, the other relevant invariants are
$x:=w, y:=v^4$ and $z:=uv$ and the quotient can be embedded in $(\bc^3,0)$.

Define $r_\gamma\in\{0,1\}$ such that $\gamma-r_\gamma=2c_\gamma$ is even. 
Notice that $\lambda=-1$.
Then, if $m=w^\alpha u^\beta v^\gamma$ is a monomial of $f^{ab}$, then 
$\beta+\gamma=2(2t_m+r_\gamma+1)$ for some $t_m\geq 0$. Set $k_0:=
\max_m(t_m-c_\gamma)$. Then 
$mv^{4k_0+2}=x^\alpha z^\beta y^{k_0-t_m+c_\gamma}$.
If $f^{ab}=\sum_ma_mm$ then the equation of the quotient is
$f':=\sum_ma_m x^\alpha z^\beta y^{k_0-t_m+c_\gamma}$.
Notice that the contribution of $u^2$ is $z^2y^{k_0}$. 
Let $d(m)$ be the $(w,u,v)$-degree $\alpha+\beta+\gamma=\alpha+4t_m+2r_\gamma
+2$ of $m$, respectively, let  $d'(m)$ be
the $(x,y,z)$-degree of $mv^{4k_0+2}$.
(In particular, $d(m)\geq 2$ with equality if and only if $m=u^2$.)
By on easy verification one gets that $d'(m)\geq 2+k_0=d'(u^2)$,
and if $d'(m)=2+k_0$ then $y^{k_0}$ divides the corresponding monomial
$x^\alpha z^\beta y^{k_0-t_m+c_\gamma}$. (In fact, the possible  monomials are
$z^2y^{k_0}, xy^{k_0+1}, zy^{k_0+1}, y^{k_0+2}$.) 
Hence the tangent cone of $f'$  is not irreducible.

\no (5) Assume that there is an $a_j$ (say $a_1$) of type $4s+2$ $(s\in\bz)$.
Then one has  the following possibilities:

\no (5.1) Set $\epsilon=diag(-1,i,i)$ (acting on $(w,u,v)$). The 
principal invariants in $\bc [w,u,v]$ are
$I=\{v^4,v^3u,v^2u^2,vu^3,u^4,wu^2,wuv,wv^2,w^2\}$. In particular, 
$embdim(\bc^3/\bz_4)=9$.

\no (5.1.1) Assume that $f^{ab}=w^2+(\deg\geq 3)$. Then 
$f^{ab}$ is an invariant, which expressed in terms of $I$
has a linear term. Hence $embdim\{f^{ab}=0\}/\bz_4=8$.

\no (5.1.2) Assume that $f^{ab}=u^2+\ldots$, in particular $\lambda=-1$.
 Then $u^4,u^3v,u^2v^2$ and $wu^2$ can be eliminated using $f^{ab}$,
and one remains with the other five principal invariants
$I_r=\{w^2,v^4,uv^3,wv^2,wuv\}$.
Write $f^{ab}$ again in the form $f^{ab}=(u+h_2+h_3)^2+(\deg\geq 5)$.
Then $uh_2$ and $h_2^2+2uh_3$ are $(-1)$-eigenfunctions, hence $h_2$
and $h_3$ are  linear combination of $w^2u$ and $w^2v$.
Hence $f^{ab}=u^2+au^2w^2+buvw^2+(\deg\geq 5)$. Multiplying such an $f^{ab}$
with any $(-1)$-eigenfunction $h$ such that $f^{ab}h$ can be expressed in terms
of the monomials $I_r$, the expression of $f^{ab}h$ in terms of these 
invariants $I_r$ will contain no linear term. Hence, $embdim
\{f^{ab}=0\}/\bz_4=5$. 

\no (5.2) Assume that $\epsilon=diag(-1,i,-i)$. The 
principal invariants  are
$I=\{v^4,vu,u^4,wu^2,wv^2,w^2\}$. 

\no (5.2.1) If $f^{ab}=w^2+(\deg\geq 3)$ then 
$f^{ab}$ is an invariant, which expressed in terms of $I$
has a linear term. Hence $embdim\{f^{ab}=0\}/\bz_4=6-1=5$.

\no (5.2.2) Assume that $f^{ab}=u^2+\ldots$, in particular $\lambda=-1$.
 Then $I_r=\{w^2,v^4,uv,wv^2\}$.
Write $f^{ab}=(u+h_2+h_3)^2+(\deg\geq 5)$.
Then $uh_2$ and $uh_3$ are $(-1)$-eigenfunctions. Analyzing the corresponding
monomial eigenfunctions, we get that 
$f^{ab}=u^2+au^2w^2+bv^2w^2+cuvw+(\deg\geq 5)$. Multiplying such an $f^{ab}$
with any $(-1)$-eigenfunction $h$ such that $f^{ab}h$ can be expressed in terms
of the monomials $I_r$, the expression of $f^{ab}h$ in terms of these 
invariants $I_r$ will contain no linear term. Hence, $embdim
\{f^{ab}=0\}/\bz_4=4$. 

(The case $f^{ab}=v^2+\cdots$ is similar.)

\no (6) Assume that all the integers $a_j$ are odd.
Then it is enough (modulo a symmetry) to consider the cases $f^{ab}=w^2+
(\deg\geq 3)$ with three different actions for $\epsilon$, namely
$diag(i,i,i)$, $diag(i,i,-i)$ and $diag(i,-i,-i)$.

In all of these cases the embedded dimension of the quotient is
$> 3$ (it is the cardinality of $I_r$, namely 9, 5, resp. 7). 
The verification is exactly the same as in (5.1.2) or (5.2.2). 
\end{proof}

\noindent 
{\bf Let us summarize} what we have: the superisolated singularity 
$(X,0)$ is clearly a Gorenstein singularity  with $b_1(M)=0$.
It has only one  ``splice equation'' (cf. \ref{UACC})
which defines $(X_w^{ab},0)$. The above fact shows that
the universal abelian cover $(X^{ab},0)$ is not in the 
equisingular deformation of $(X_w^{ab},0)$. Therefore, {\em even if the 
construction of the ``splicing equations'' is not obstructed
(cf. \ref{4.2} and \ref{4.3}), in general, the $UACC$
\cite{NWnew} is not valid.}

We can go even further: the resolution graphs of $(X^{ab},0)$ and
$(X_w^{ab},0)$ are the same, hence these two singularities
 have the same topological types. Their links are rational homology 
spheres (with first homology $\bz_7$). Since the common resolution graph
is star-shaped, and $(X_w^{ab},0)$ is weighted homogeneous, $(X^{ab},0)$
can be equisingularly deformed into $(X_w^{ab},0)$ if and only if
their geometric genus are the same (cf. with the already mentioned 
result of Neumann and Wahl \cite{NW} (3.3)). Since this is not the case
(by the above Fact A), one gets that $p_g(X^{ab},0)\not=
p_g(X_w^{ab},0)$. In particular, {\em we constructed two Gorenstein
singularities (one of them is even a hypersurface Brieskorn singularity)
with the same rational homology sphere link, but with different geometric
genus.} This provides counterexample for both $SWC$ and $GGC$. 

What is even more striking in the above counterexample, is the fact
that the corresponding graphs are star-shaped (and one of the singularity
is weighted homogeneous), cf. with the last paragraph of \ref{GGC}. 

Recall that $p_g(X_w^{ab},0)=10$.   Notice also 
that for {\em any} normal surface singularity with the same resolution 
graph as $\Gamma^{ab}$, by (9.6) of \cite{NOSZ} one has $p_g\leq 10$. 
In particular, $p_g(X^{ab},0)<10$. 

\subsection{Counterexample: The case $C_2$ with  $d=5$ and multiplicity 
sequence $[2_6]$, and  its universal abelian cover.}\label{4.5} \
We start with $f=f_5+z^6$ where
$f_5=z(yz-x^2)^2-2xy^2(yz-x^2)+y^5$.
The curve $C$ is irreducible with unique 
singularity at $[0:0:1]$ (of type $A_{12}$).
The resolution graph $\Gamma$ of the superisolated singularity $(X,0)$
is

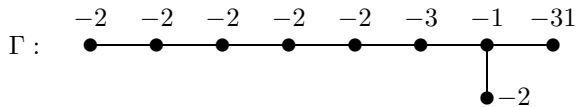
\begin{picture}(300,45)(-30,0)
\put(25,25){\circle*{5}}
\put(50,25){\circle*{5}}
\put(75,25){\circle*{5}}
\put(25,35){\makebox(0,0){$-2$}}
\put(50,35){\makebox(0,0){$-2$}}
\put(75,35){\makebox(0,0){$-2$}}
\put(100,25){\circle*{5}}
\put(125,25){\circle*{5}}
\put(150,25){\circle*{5}}
\put(175,25){\circle*{5}}
\put(200,25){\circle*{5}}
\put(175,5){\circle*{5}}
\put(25,25){\line(1,0){175}}
\put(175,25){\line(0,-1){20}}
\put(0,25){\makebox(0,0){$\Gamma:$}}
\put(100,35){\makebox(0,0){$-2$}}
\put(125,35){\makebox(0,0){$-2$}}
\put(150,35){\makebox(0,0){$-3$}}
\put(175,35){\makebox(0,0){$-1$}}
\put(200,35){\makebox(0,0){$-31$}}
\put(185,5){\makebox(0,0){$-2$}}
\end{picture}

Since the graph is star-shaped, the same resolution graph
can be realized by a weighted homogeneous singularity $(X_w,0)$ as well.
In fact, it is much easier to determine the universal abelian cover
$(X_w^{ab},0)$ of $(X_w,0)$. By \cite{Neu}, it is the Brieskorn
hypersurface singularity $\{u^{13}+v^{31}+w^2=0\}$ (and this agrees with the
``splice equation'' provided by $\Gamma$).
The corresponding resolution graph $\Gamma^{ab}$ (of both $(X^{ab},0)$
and $(X_w^{ab},0)$) is 

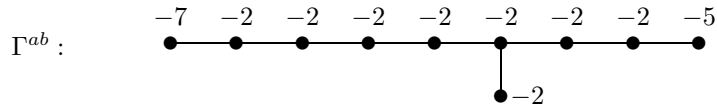
\begin{picture}(300,45)(0,0)
\put(0,25){\makebox(0,0){$\Gamma^{ab}:$}}
\put(50,25){\circle*{5}}
\put(75,25){\circle*{5}}
\put(100,25){\circle*{5}}
\put(125,25){\circle*{5}}
\put(150,25){\circle*{5}}
\put(175,25){\circle*{5}}
\put(200,25){\circle*{5}}
\put(225,25){\circle*{5}}
\put(250,25){\circle*{5}}
\put(175,5){\circle*{5}}
\put(50,25){\line(1,0){200}}
\put(175,25){\line(0,-1){20}}
\put(50,35){\makebox(0,0){$-7$}}
\put(75,35){\makebox(0,0){$-2$}}
\put(100,35){\makebox(0,0){$-2$}}
\put(125,35){\makebox(0,0){$-2$}}
\put(150,35){\makebox(0,0){$-2$}}
\put(175,35){\makebox(0,0){$-2$}}
\put(200,35){\makebox(0,0){$-2$}}
\put(225,35){\makebox(0,0){$-2$}}
\put(250,35){\makebox(0,0){$-5$}}
\put(185,5){\makebox(0,0){$-2$}}
\end{picture}

\no which defines the 
{\em integral homology sphere}  $\Sigma(13,31,2)$.
%splice diagram 

%\begin{picture}(300,45)(0,0)
%\put(150,25){\circle*{5}}
%\put(200,25){\circle*{5}}
%\put(175,25){\circle*{5}}
%\put(175,5){\circle*{5}}
%\put(150,25){\line(1,0){50}}
%\put(175,25){\line(0,-1){20}}
%\put(165,35){\makebox(0,0){$13$}}
%\put(185,35){\makebox(0,0){$31$}}
%\put(180,15){\makebox(0,0){$2$}}
%\end{picture}

The action of $H=\bz_5$ on $(X_w^{ab},0)$
is $(u,v,w)\mapsto (\zeta^4u, \zeta^2v,\zeta w)$,
where $\zeta$ denotes a 5-root of unity. This action  has a lot of principal
invariants, but one can eliminate those ones which are multiples of
$w^2$ using the equation 
$u^{13}+v^{31}+w^2$. Therefore, we have to consider only the following 
ones: $a:=u^5,\ b:=v^5,\ c:=u^2v,\ d:=uv^3,\ e:=uw$ and $f:=wv^2$.
If one wants to get the equations of $(X_w,0)$ in $\bc^6$ (in variables
$a,\cdots,f$), one has to eliminate from the equations
$u^{13}+v^{31}+w^2,u^5-a,v^5-b,u^2v-c,uv^3-d,uw-e,wv^2-f$
the variables $(u,v,w)$. This can be done by {\sc Singular} \cite{GPS01}, and we get
the following set of equations for $X_w$ in $\bc^6$:
$$X_w=\left\{ 
\begin{array}{ll}
ab-c^2d=0 & \\
bc-d^2=0 & \\
ad-c^3=0 & \\
be-df=0 & \\
de-cf=0 & \\
af-c^2e=0 & \\
e^2+a^3+b^6c=0 & \\
ef+a^2c^2+b^6d=0 & \\
f^2+ac^4+b^7=0 & \end{array}\right.$$
In fact, these equations can also be obtained without {\sc Singular}:
the first six equations are the principal relations connecting
the principal invariants $a, \ldots, f$, while the last three equations 
are obtained (see the recipe in the proof of Fact B, step (3))
by multiplying the $\zeta^2$-eigenfunction $u^{13}+v^{31}+w^2$
by the $\zeta^3$-eigenfunctions $u^2,uv^2,v^4$.

As a curiosity, separating the first  six equations one gets 
that $(X_w,0)$ is a subgerm of the determinantal singularity
defined by the $(2\times 2)$-minors of 
$$\left(\begin{array}{llll} b&d&f&c^2\\d&c&e&a \end{array}\right).$$

The weights of the variables $(a,\ldots,f)$ are (62,26,30,28,93,91).

\vspace{2mm}

Notice also that $(X_w,0)$ is Gorenstein, but it is not a
complete intersection. Moreover, 
{\em the two singularities $(X,0)$ and $(X_w,0)$ 
have the same topological types (the same graphs $\Gamma$),
but their embedded dimensions are not the same}: they are 
3 and 6 respectively. It is even more surprising that {\em
their multiplicities are also different}: $mult(X,0)=5$ and
$mult(X_w,0)=6$ (the second computed by {\sc Singular} \cite{GPS01}). 

\vspace{2mm}

On the other hand, their geometric genera are the same: 
$p_g(X,0)=10$ by the formula of (4), $p_g(X_w,0)=10$ by
Pinkham's formula \cite{Pi1}. In particular, using again
\cite{NW} (3.3), $(X,0)$ is in the equisingular deformation 
of $(X_w,0)$. 

This deformation can be described as follows.
(Again, the weight of $\lambda $ is $-3$.) 
The authors are grateful
to J. Stevens for his help in finding these deformation.
$$X(\lambda)=\left\{ 
\begin{array}{ll}
ab-c^2d=\lambda f & \\
bc-d^2=\lambda^2 a & \\
ad-c^3=\lambda e & \\
be-df=-\lambda ac^2 & \\
de-cf=-\lambda a^2 & \\
af-c^2e=-\lambda b^6 & \\
e^2+a^3+b^6c=0 & \\
ef+a^2c^2+b^6d=0 & \\
f^2+ac^4+b^7=0 & \end{array}\right.$$
In order to understand the deformation, consider the equation (for
$\lambda\not=0$):
$$E:=\lambda^{-2}( a^2b-2ac^2d+c^5)+b^6.$$
Notice that for $\lambda\not=0$, using the first three equations
one can eliminate the variables $a,e,f$.
The last four  equations transform into $E\lambda$,
$Ec$, $Ed$ and $Eb$ (where in $E$
we substitute $a$).
Hence their vanishing is equivalent with the vanishing of $E$. 
The forth and fifth equations are automatically satisfied.
Hence, for $\lambda\not=0$,  the system of equation is equivalent 
with a hypersurface singularity in variables $(b,c,d)$ given by
$E=0$ with the substitution $a=\lambda^{-2}(bc-d^2)$.
Taking $\lambda =1$, $b=z$, $c=y$ and $d=x$, one gets exactly
the superisolated singularity $f=f_5+z^6$.

\vspace{2mm}

On the other hand, similarly as in the case of \ref{4.4},
{\em there is no equisingular deformation 
at the level of universal abelian covers.}
Both $(X^{ab},0)$ and $(X^{ab}_w,0)$ have the same graph $\Gamma^{ab}$
-- which is a {\em unimodular star-shaped graph}, but $(X^{ab},0)$
is not in the equisingular deformation of $(X_w^{ab},0)$.
In particular (by the same argument as in \ref{4.4}),
$p_g(X^{ab},0)< p_g(X_w^{ab},0)$. In particular, {\em all the conjectures
UACC, SWC  and GGC  fail.}
(For the first case notice that the ``splice equation'' of $(X,0)$ 
is exactly the equation of $(X^{ab}_w,0)$.) 

This example shows (cf. with \ref{4.4})
that even with the assumption $H=0$ counterexamples for these 
conjectures exist.

\vspace{2mm}

The non-existence of the deformation follows by a similar statement as in 
the case of \ref{4.4}:
{\em Assume that $\bz_5$ acts
freely in codimension 1 on a hypersurface germ which in some coordinates
has the form $w^2+(\deg \geq 6)$. Then if 
the quotient is a hypersurface with multiplicity greater than $2$,
then the tangent cone  of the quotient is reducible. }

This has a completely similar proof as the similar statement in
\ref{4.4}, and we will not give it here. 

\section{Integral homology sphere links}

\subsection{} Recall that the first homology of the 
link of a hypersurface superisolated
singularity $f=f_d+l^{d+1}$ is $\bz_d$ ($d\geq 2$); in particular,
it is never trivial. Nevertheless, 
we would like to emphasize that even with 
integral homology sphere links, counterexamples exist (although it is 
harder to find them); this additional requirement
does not change the picture. But, in order to find such examples, 
we have to enlarge our
family. We exemplify here two possibilities.

\subsection{The universal abelian cover revisited.} \
In the first case we consider the universal abelian cover $(X^{ab},0)$
of a hypersurface superisolated singularity. Notice that, in general,
it is hard to give the equations (or identify the analytic structure)
of $(X^{ab},0)$. But its topological type can be described completely.
Recall that the minimal resolution $\tilde{X}$ of $(X,0)$ contains 
only one exceptional divisor $C$ with (plane curve) singularities $(C,p_i)$
and self-intersection $-d$. If one considers the $\bz_d$-cyclic 
cover $q:\tilde{X}^{ab}\to \tilde{X}$
of $\tilde{X}$, branched along $C$, one gets a partial resolution 
of $(X^{ab},0)$. In general $\tilde{X}^{ab}$ is not smooth, its singularities
are the $d$-suspensions of the plane curve singularities $(C,p_i)$
(in other words, if the local equation of $(C,p_i)$ is $g_i(u,v)=0$,
then $Sing\tilde{X}^{ab}=q^{-1}(\cup_ip_i)$, $q^{-1}(p_i)$ contains only one 
point, and $(\tilde{X}^{ab},q^{-1}(p_i))$ is a 
hypersurface singularity of type
$g_i(u,v)+w^d=0$).
Moreover, the self-intersection of $\tilde{C}:=q^{-1}(C)$ is $-1$.
(Indeed, $d\tilde{C}\cdot \tilde{C}= q^*C\cdot \tilde{C}=C\cdot q_*\tilde{C}
=C^2=-d$.) 
In particular, the minimal good resolution graph of $(X^{ab},0)$
can be obtained in a similar way as the graphs of hypersurface
suspension singularities (if one replaces the embedded resolution
graphs of the plane curve singularities with the graphs of their 
$d$-suspensions, and the self-intersection $-d$ with $-1$). 

This construction also shows that the link of $(X^{ab},0)$ is an
integral homology sphere if and only if all the links of the 
$d$ suspension singularities $g_i(u,v)+w^d$ are integral homology
spheres. This fact can be realized, as it is shown by the example
\ref{4.5}. But even with $N\geq 3$ one can find many examples.

Take for example $C_8$ in the table \ref{4.1} with $d=5$ and $N=4$.
Then the local equations of the plane curve singularities are $u^7+v^2$ 
and three times $u^3+v^2$. Hence the minimal good resolution graph
of $(X^{ab},0)$ is unimodular, and has the following form
(where all the undecorated curves are $-2$-curves):

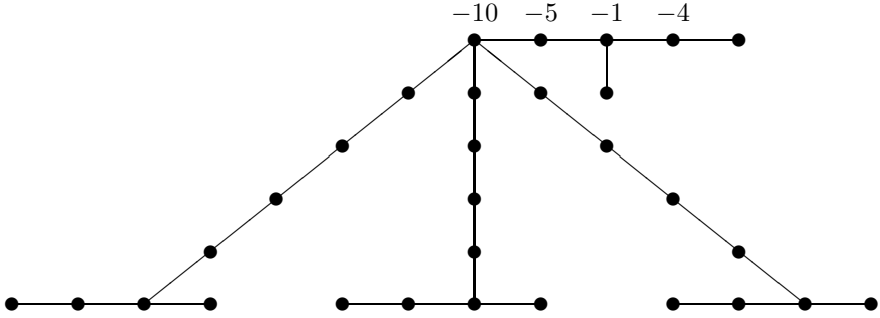
\begin{picture}(300,120)(20,0)
\put(25,5){\circle*{5}}
\put(50,5){\circle*{5}}
\put(75,5){\circle*{5}}
\put(100,5){\circle*{5}}
\put(150,5){\circle*{5}}
\put(175,5){\circle*{5}}
\put(200,5){\circle*{5}}
\put(225,5){\circle*{5}}
\put(275,5){\circle*{5}}
\put(300,5){\circle*{5}}
\put(325,5){\circle*{5}}
\put(350,5){\circle*{5}}
\put(100,25){\circle*{5}}
\put(200,25){\circle*{5}}
\put(300,25){\circle*{5}}
\put(125,45){\circle*{5}}
\put(200,45){\circle*{5}}
\put(275,45){\circle*{5}}
\put(150,65){\circle*{5}}
\put(200,65){\circle*{5}}
\put(250,65){\circle*{5}}
\put(175,85){\circle*{5}}
\put(200,85){\circle*{5}}
\put(225,85){\circle*{5}}
\put(250,85){\circle*{5}}
\put(200,105){\circle*{5}}
\put(225,105){\circle*{5}}
\put(250,105){\circle*{5}}
\put(275,105){\circle*{5}}
\put(300,105){\circle*{5}}
\put(25,5){\line(1,0){75}}
\put(150,5){\line(1,0){75}}
\put(275,5){\line(1,0){75}}
\put(75,5){\line(5,4){125}}
\put(200,5){\line(0,1){100}}
\put(325,5){\line(-5,4){125}}
\put(200,105){\line(1,0){100}}
\put(250,105){\line(0,-1){20}}
\put(200,115){\makebox(0,0){$-10$}}
\put(225,115){\makebox(0,0){$-5$}}
\put(250,115){\makebox(0,0){$-1$}}
\put(275,115){\makebox(0,0){$-4$}}
\end{picture}

This shows that, even if we deal with integral homology sphere links,
in general, the {\em semigroup conjecture fails} (cf. \ref{UACC}). 

\subsection{Non-hypersurfaces.} Another 
way to extend our class of examples is to
consider all the singularities (not only the hypersurfaces)
which have the property that one of their resolution graphs has 
a ``central vertex'', and all the graph-components of 
the complement of the central vertex  are embedded resolution graphs of 
plane curve singularities. 

Probably the simplest (non-hypersurface)  example is the following
complete intersection in $(\bc^4,0)$, given by the equations
$$(X,0)=\{x^2=u^3+v^2y,\ y^2=v^3+u^2x\}.$$
Its resolution graph is unimodular and has the form:

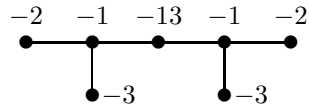
\begin{picture}(300,45)(0,0)
\put(125,25){\circle*{5}}
\put(150,25){\circle*{5}}
\put(175,25){\circle*{5}}
\put(200,25){\circle*{5}}
\put(225,25){\circle*{5}}
\put(150,5){\circle*{5}}
\put(200,5){\circle*{5}}
\put(125,25){\line(1,0){100}}
\put(150,25){\line(0,-1){20}}
\put(200,25){\line(0,-1){20}}
\put(125,35){\makebox(0,0){$-2$}}
\put(150,35){\makebox(0,0){$-1$}}
\put(175,35){\makebox(0,0){$-13$}}
\put(200,35){\makebox(0,0){$-1$}}
\put(225,35){\makebox(0,0){$-2$}}
\put(160,5){\makebox(0,0){$-3$}}
\put(210,5){\makebox(0,0){$-3$}}
\end{picture}

This example appears in \cite{NW} as a ``positive'' example
satisfying the Casson invariant conjecture.

In the spirit of the this section, we present  one of its ``negative''
properties:  {\em its minimal (Artin) cycle does not agree with its
maximal cycle (i.e. the minimal cycle cannot be cut out by
a holomorphic function-germ).}
Notice that, in general, if one wishes
a topological characterization of the multiplicity, the first test is
exactly the identity of the minimal and maximal cycles.

In order to see that in this case they are not the same, 
notice two facts. First, the 
strict transforms of the four coordinate functions are supported by the 
four leaves (degree  one vertices). Second,
the intersection of the minimal cycle with $C$ (the $-13$-curve)
is $-1$, and with 
all the other irreducible exceptional divisors is zero. 
In particular, analyzing the graph (e.g. 
the corresponding linking numbers), one gets that if the divisor of 
a holomorphic germ  $f$ would be the sum of the minimal cycle and the strict
transform, then the local intersection multiplicity $i_{(X,0)}(f,z)$
would be 2 for one (in fact, for two) of the coordinate functions $z$.
This would imply that the multiplicity of $(X,0)$ is 2 (or less), 
in particular $(X,0)$ would be 
a hypersurface singularity. But this is not the case. 
(Nevertheless, the topological type of $(X,0)$ supports 
at least one analytic structure
 for which the maximal cycle is the minimal cycle.)


\begin{thebibliography}{99}

%\bibitem{Casson} Akbulut, Selman and McCarthy, John D.: Casson's invariant
%for oriented homology 3--spheres, an exposition, {\em 
%Mathematical Notes} {\bf 36},
%Princeton University Press, Princeton, 1990.

%\bibitem{Arndt} Arndt, J.: Verselle Deformationen zyklischer 
%Quotientensingularit\"aten, Diss. Hamburg, 1988. 

%\bibitem{AGV} Arnold, V.I and  Gusein-Zade, S.M. and Varchenko, A.N.: 
%{\em Singularities of Differentiable maps}, Volume 1 and 2, Monographs Math.,
%{\bf 82-83}, Birkh\"auser, Boston, 1988.

\bibitem{Artal} Artal  Bartolo, E.: Forme de Jordan de la monodromie  des 
singularit\'es superisol\'ees de surfaces,
{\it Mem. Amer. Math. Soc.}, {\bf 525} (1994). 


\bibitem{aclm} Artal Bartolo, E.; Cassou-Nogu\`es, Pi.; Luengo, I.; 
Melle Hern\'andez, A.:
 Monodromy conjecture for some surface singularities,  
{\em Ann. Sci. \'Ecole Norm. Sup. (4)}  {\bf 35}  (2002),  no. 4, 605--640.

%\bibitem{AB} Artal--Bartolo, Enrique: Forme de Seifert des singularit\'es de
% surface, {\em C. R. Acad. Sci. Paris}, {\bf t. 313}, S\'erie I, 689-692,
%1991.

\bibitem{Artin62} Artin, M.: 
Some numerical criteria for contractibility of curves on algebraic surfaces,
{\em  Amer. J. of Math.}, {\bf 84}, 485-496, 1962.

\bibitem{Artin66} Artin, M.: 
On isolated rational singularities of surfaces,
{\em Amer. J. of Math.}, {\bf 88}, 129-136, 1966.

%\bibitem{Art} Artin, M.: Algebraic construction of Brieskorn's 
%resolutions, {\em J. of Algebra} {\bf 29} (1974), 330-348.

%\bibitem{A} Ashikaga, Tadashi: The Signature of the Milnor fiber of
%Complex Surface Singularities on Cyclic Coverings, preprint, 1995.
%
%\bibitem{At} Atiyah, Michael F.: New invariants of 3-- and 4--manifolds,
%{\em The Mathematical  Heritage of Herman Weyl}, {\em Proc. Symp. Pure
%Math.} {\bf 48}, 285-300, 1988.
%
%\bibitem{APS} Atiyah, M.F., Patodi, V.K. and Singer, I.M.: Spectral asymmetry
%and Riemannian geometry, I,II,III, {\em Math. Proc. Cambridge Philos. Soc.},
%{\bf 77} 53-69, 1975; {\bf 78} 405-432, 1975; {\bf 79} 71-99, 1976.
%
%\bibitem{B1} B\u adescu, Lucian: 
%Applications of the Grothendieck duality theory
%to the study of normal isolated singularities,
%{\em Rev. Roum. Math. Pures et
%Appl.}, Tome {\bf XXIV}, No 5, 673-689, 1979.
%
%\bibitem{B2} B\u adescu, Lucian: Dualising divisors of two--dimensional
%singularities,
%{\em Rev. Roum. Math. Pures et Appl.}, Tome {\bf XXV}, No 5, 695-707, 1980.

%\bibitem{BK} Behnke, K. and Kn\"orrer, H.: On infinitesimal deformations of 
%rational surface singularities, {\em Comp. Math.}, {\bf 61} 1987), 103-127.

%\bibitem{BoNi} Boyer, S. and Nicas, A.:
%Varieties of group representations
%and Casson's invariant for rational homology 3--spheres, 
%{\em Trans. AMS} {\bf 322} (2), 507-522, 1990.

%\bibitem{BS1} Brian\c{c}on J. and Speder, J.: La trivialit\'e topologique
%n'implique pas les conditions de Whitney, {\em C. R. Acad. Sc., Paris}, 
%{\bf 280} (1975), 365-367.

%\bibitem{BS} Brian\c{c}on J. and Speder, J.: Les condition de Whitney 
%impliquent  ``$\mu^*$ constant'', {\em Ann. Inst. Fourier (Grenoble)}, 
%{\bf 26} (1976), 153-164.

%\bibitem{Brtop} Brieskorn, E.: Examples of singular normal complex spaces 
%which are topological manifolds, {\em Proc. Nat. Acad. Sci.}.
%{\bf 55} (1966), 1395-1397.

%\bibitem{Brtop2} Brieskorn, E.: Beispiele zur Differentialtoplogie von
%Singularit\"aten, {\em Invent. Math.}, {\bf 2} (1966), 1-14.

%\bibitem{Brsing} Brieskorn, E,: Singular elements in semi-simple algebraic
%groups, {\em Proc. Int. Con. Math. Nice}, {\bf 2} (1971), 279-284.

%\bibitem{Chen1} Chen, W:  Casson invariant and Seiberg-Witten gauge 
%theory, {\em Turkish J. Math.}, {\bf 21}(1997), 61-81.

%\bibitem{CG} Christophersen, J. A. and Gustavsen, T. S.: On infinitesimal
%deformations and obstructions for rational surface singularities,
%{\em J. Algebraic Geometry}, {\bf 10} (2001) (1), 179-198.

%\bibitem{[44]} Chunsheng Ban,  McEwan, Lee  and N\'emethi, Andr\'as:
%The embedded resolution of $f(x,y)+z^2:
%({\bf C}^3,0)\to ({\bf C},0)$,
%{\em Studia Scientiarum Math. Hungarica}, {\bf 38}, 2001, 51-71.

\bibitem{CoS} Collin, O. and Saveliev, N.: A geometric proof of the
Fintushel-Stern formula, {\em Adv. in Math.}, {147}(1999), 304-314 .

%\bibitem{CoS2} Collin, O. and Saveliev, N.: Equivariant Casson invariants
%via gauge theory, {\em J. reine Angew. Math.}, {\bf 541} (2001), 143-169.

\bibitem{CoS3} Collin, O. and Saveliev, N.:
 Equivariant Casson invariant for knots and the
Neumann-Wahl formula, {\em Osaka J. Math.}, {37} (2000), 57-71.

%\bibitem{Di1} Dixon, D. : Dissertation, Purdue University, 1977.

%\bibitem{Di2} Dixon, D.J.: The fundamental divisor of normal double
%points of surfaces,
%{\em Pacific J. Math.}, {\bf 80}, no. 1, 105-115, 1979.

%\bibitem{Dolg} Dolgachev, I.V.: Automorphic forms and weighted homogeneous
%singularities, {\em Funkt. Anal. Jego. Prilozh.}, {\bf 9}(1975), 67-68 .
%English translation in {\em Funct. Anal. Appl.}, {\bf 9}(1975), 149-151.

%\bibitem{Du} Durfee, A.: The Signature of Smoothings of Complex
%Surface Singularities, {\em Math. Ann.}, {\bf 232}, 85-98, 1978.

%\bibitem{Ebe} Ebeling, W.: Poincar\'{e} series and monodromy of a 
%two-dimensional quasihomogeneous hypersurface singularity, 
%{\sf math.AG/0109210}, {\em Manuscripta Math.}, to appear.

%\bibitem{EN} Eisenbud, D. and Neumann, W.:
%{\em Three-Dimensional Link Theory
%and Invariants of Plane Curve Singularities}, Ann. of Math. Studies
%{\bf 110}, Princeton University Press, 1985.

%\bibitem{Er} Elkik, R.: Singularit\'es rationelles et D\'eformations,
%{\em Inv. Math.}, {\bf 47} (1978), 139-147. 

\bibitem{Fenske} Fenske, T.: Rational 1- and 2-cuspidal plane curves, {\em Beitr\"age Algebra Geom.} {\bf 40}  (1999), no. 2, 309--329. 

\bibitem{4authors} Fern\'andez de Bobadilla, J.; Luengo-Valesco, I.;
Melle-Hern\'andez, A. and N\'emethi, A.: manuscript in preparation. 

\bibitem{FS} Fintushel, R. and Stern, R.J.: Instanton homology of
Seifert fibered homology 3--spheres, {\em Proc. London Math. Soc.}, (3)
{\bf 61}, 109-137, 1991.

\bibitem{fujita} Fujita, G.:  A splicing  formula for the Casson-Walker's
invariant, {\em  Math. Annalen}, {\bf 296}, 327-338 (1993).

%\bibitem{NG1}
%Garc\'\i a L\'opez, Ricardo  and N\'emethi, Andr\'as:
% On the monodromy at infinity of a polynomial map,
%{\em Compositio  Math.}, {\bf 100}, 205-231, 1996.
%Appendix by Garc\'\i a  L\'opez, Ricardo  and Steenbrink, Joseph.

%\bibitem{NG2} Garc\'\i a L\'opez, Ricardo and N\'emethi, Andr\'as:
%On the monodromy at infinity of a polynomial map, II.
%{\em Compositio Math. }, {\bf 115}, 1-20, 1999.

%\bibitem{NG3} Garc\'\i a L\'opez, Ricardo and N\'emethi, Andr\'as:
%Hodge numbers attached to a polynomial map,
%to appear in {\em Ann. Inst. Fourier, Grenoble}. \
%http://www.mi.aau.dk/\~\,esn/titles\_1997.html

%\bibitem{GM} M. Goresky, R. MacPherson: Intersection Homology Theory,
%{\em Topology}, {\bf 19},  135-162, 1980. 
%

%\bibitem{Grauert} Grauert, H.: \"Uber Modifikationen  und exceptinelle
%analytische Mengen, {\em Math. Annalen}, {\bf 146} (1962), 331-368.

%\bibitem{GrauDef} Grauert, H.: \"Uber die Deformationen Isolierten 
%Singularit\"aten Analytischer Mengen, {\em Inv. Math.}, {\bf 15} (1972), 
%171-198. 

%\bibitem{Greenberg} Greenberg, Leon: Homomorphisms of triangle groups into
%$PSL(2,\bc)$, {\em Riemannian surfaces and related topics: Proc. of the 1978
%Stony Brook Conference}, Annals of Math. Studies {\bf 97}, Princeton
%University Press 1981, 167-181.

%\bibitem{Greuel}  Greuel, G.-M.:  Constant Milnor number implies constant
%multiplicity for quasi-homogeneous singularities, {\em Manuscripta
%Math.}, {\bf 56} (1986), 159-166.

%\bibitem{GH} Greuel G.-M. and  Hamm, H.A.:  Invarianten 
%quasihomogener vollst\"andiger Durchschnitte, {\em Invetiones math.}, {\bf 
%49} (1978), 67-86.

%\bibitem{GSt} Greuel,  G.-M. and  Steenbrink, J.: On the topology of 
%smoothable singularities, {\em Proc. of Symp. in Pure Math.}, 
%{\bf 40} Part 1 (1983), 535-545.

%\bibitem{CDG} Gusein-Zade, S.M., Delgado, F. and Campillo, A.:
%On the monodromy of a plane curve singularity  and the Poincar\'e 
%series of the ring of functions on the curve, 
%{\em Functional Analysis and its Applications}, {\bf 33}(1)  (1999), 56-67.

%\bibitem{Hamm} Hamm. H.A.: Exotische Sph\"aren als Umgebungsr\"nder in
%speziellen komplexen R\"aumen, {\em Math. Ann.}, {\bf 197} (1972), 44-56.

%\bibitem{HKK} Harer, John, Kas, Arnold and Kirby, Robion:
%Handlebody decomposition of complex surfaces, {\em Mem. Amer. Math. Soc.},
%{\bf 62} (1986) no 350.

%\bibitem{Report} Hauser, H.  and Randell, R.:
%Report on the problem session, {\em Contemporary Mathematics}, Vol. {\bf 90},
%119-134, 1989. 

%\bibitem{Ishii} Ishii, S:  The invariant $-K^2$ and continued fractions
%for 2-dimensional cyclic quotient singularities, preprint.

%\bibitem{JN} Jankins, M. and  Neumann, W.D.: {\sl Lectures on Seifert
%Manifolds}, Brandeis Lecture Notes, 1983.

%\bibitem{JongStraten1} de Jong, T. and van Straten, D.: On the base space 
%of a semi-universal deformation of rational quadruple points, {\em 
%Annals of Math.}, {\bf 134} (2) (1991), 653-678.
 
%\bibitem{JongStraten2} de Jong, T. and van Straten, D.: On the 
%deformation theory  of rational surface singularities with reduced 
%fundamental cycle, {\em J. Alg. Geom.}. {\bf 3} (1994), 117-172.

%\bibitem{JongStraten} de Jong, T. and van Straten, D.: Deformation theory
%of sandwiched singularities, {\em Duke Math. J.}, {\bf 95} (3) (1998), 
%451-522.

%\bibitem{K1} Karras, U.: On Pencils of Curves and Deformations of Minimally
%Elliptic Singularities, {\em Math. Ann.}, {\bf 247}, 43-65, 1980.

%\bibitem{K2} Karras, U.: Methoden zur Berechnung algrebraischer 
%Invarianten und zur
%Konstruktion von Deformationen normaler Fl\"achensingularit\"aten,
%Habilitationsschrift, Dortmund 1981.

%\bibitem{Kaw} Kawamata, Yujiro: Crepant blowing--up of 3--dimenasional
%canonical singularities and its application to degeneration of surfaces,
%{\em Annals of Math.}, {\bf 127}, 93-163, 1988. 

%\bibitem{Kn} Kn\"oller, F.W.: 2--dimensionale Singularit\"aten und
%Differentialformen, {\em Math. Ann.}, {\em 206} (1973), 205-213.

%\bibitem{KrSt} Kreck, M. and  Stolz, S.: Nonconnected moduli spaces
%of positive sectional curvature metrics, {\em J. of the AMS.}, 
%{\bf 6}(1993), 825-850.

%\bibitem{KSB} Koll\'ar, J. and Shepherd-Barron, N.I.: Threefolds and
%deformations of surface singularities, {\em Invent. math.}, {\bf 91}
%(1988), 299-338.

%\bibitem{flips} Koll\'ar, J.: Flips, flops, minimal models, etc.,
%{\em Surveys in Diff. Geom.}, {\bf 1} (1991), 113-199.

%\bibitem{Kollar} Koll\'ar, J\'anos: Shafarevich Maps and Automorphic Forms, 
%Princeton University Press, Princeton, 1995.

%\bibitem{MoKo} Koll\'ar, J\'anos and Mori, Shigefumi (with the 
%collaboration of C. H. Clemens and A. Corti): Birational Geometry
%of Algebraic Varieties, Cambridge University Press, {\bf 134}, 1998.

%\bibitem{Lauferbook} Laufer, H.B.: Normal two--dimensional singularities.
%{\em Annals of Math. Studies}, {\bf 71}, Princeton University Press, 1971.

\bibitem{Laufer72} Laufer, H.B.: On rational singularities,
{\em Amer. J. of Math.}, {\bf 94}, 597-608, 1972.

\bibitem{Laufer73} Laufer, H.B.: Taut two--dimensional singularities,
{\em Math. Ann.}, {\bf 205}, 131-164,  1973.

\bibitem{Laufer77} Laufer, H.B.: On minimally elliptic singularities,
{\em Amer. J. of Math.}, {\bf 99}, 1257-1295, 1977.

\bibitem{Laufer77b} Laufer, H.B.: On $\mu$ for surface singularities,
{\em Proceedings of Symposia in Pure Math.}, {\bf 30}, 45-49,  1977.

%\bibitem{Lauferdouble} Laufer, H.B.: On normal two--dimensional double
%point singularities, {\em  Israel Journal of Mathematics}, Vol. {\bf 31}, 
%Nos 3-4, 315-334,  1978.

%\bibitem{LaW} Laufer, H.B.: Weak simultaneous resolution for deformations 
%of Gorenstein surface singularities, {\em Proc. of Symp. in Pure Math.}, 
%{\bf 40}, Part 2 (1983), 1-29. 

%\bibitem{LaS} Laufer, H.B.: Strong Simultaneous Resolution for 
%Surface Singularities, {\em Adv. Studies in Pure Math.}, {\bf 8} (1986),
%207-214. {\em Complex Analytic Singularities}. 

%\bibitem{Laumult} Laufer, H.B.: 
%The multiplicity of isolated two--dimensional
%hypersurface singularities, {\em Transactions of the AMS}, {\bf 302} 
%Number 2, 489-496, 1987.

%\bibitem{le} L\^e D\~ung Tr\'ang: Topologie des singularit\'es  des
%hypersurfaces complexes, {\em Ast\'erisque}, {\bf 7-8} (1973), 171-182.

\bibitem{Lescop} Lescop, C.: Global Surgery Formula for the Casson-Walker
Invariant, {\em Annals of Math. Studies}, vol. {\bf 140}, Princeton University
Press, 1996.

%\bibitem{Lim} Lim, Y:  Seiberg-Witten invariants for 3-manifolds in the
%case $b_1=0$ or $1$, {\em Pacific J. of Math.}, {195}(2000), 179-204.

%\bibitem{L1} Looijenga, E.: The smoothing components of a triangle 
%singularity. I, {\em Proc. of Symp. in Pure Math.}, {\bf 40 } Part 2, (1983),
%173-183.

%\bibitem{L2} Looijenga, E.: Riemann-Roch and smoothing of singularities,
%{\em Topology} {\bf 25} (3) (1986), 293-302.

%\bibitem{LW}  Looijenga, E. and Wahl, J.:  Quadratic functions and smoothing 
%surface singularities, {\em Topology}, {\bf 25} (1986), 261-291.

\bibitem{Ignacio} Luengo, I.: The $\mu$-constant stratum is not
smooth, {\em Invent. Math.}, {\bf 90} (1), 139-152, 1987.

%\bibitem{MW}  Marcolli, M. and  and  Wang, B.L.: Seiberg-Witten invariant and 
%the Casson-Walker invariant for rational homology 3-spheres,
%{\sf math.DG/0101127}, {\em Geometri{\ae} Dedicata}, to appear.

\bibitem{Melle}  Melle-Hern\'andez, A.: Milnor numbers for surface 
singularities,  {\em Israel J. Math.}, {\bf  115}  (2000), 29--50.

%\bibitem{[56]} Mendris, R.  and N\'emethi, A.:
%  The link of $\{f(x,y)+z^n=0\}$ and Zariski's Conjecture,   submitted.
%{\sf math.AG/0207212.}

%\bibitem{MBook} Milnor, J.: Singular points of complex hypersurfaces,
%{\em Annals of Math. Studies} {\bf 61}, Princeton University Press 1968.

%\bibitem{Mo} Mordell, L.J.: Lattice points in a tetrahedron and generalized
%Dedekind sums, {\em J. Indian Math.}, {\bf 15}, 41-46,  1951.

%\bibitem{Mumford} Mumford, D.: The topology of normal singularities
%of an algebraic surface and criterion for simplicity, {\em IHES
%Publ. Math.} {\bf 9}, 5-22, 1961. 

\bibitem{Namba} Namba, M.: {\em Geometry of projective algebraic curves},
 Monographs and Textbooks in Pure and Applied Mathematics, 
{\bf 88} Marcel Dekker, Inc., New York, (1984).

%\bibitem{N114} N\'emethi, Andr\'as:
%Algebraic torsion, zeta function and Dirichlet series for graph
%       links in homology 3-spheres,
%        {\em Ergodic Theory and Dynamical Systems,} {\bf 13}, 131-142, 1993.
%
%\bibitem{N116} N\'emethi, Andr\'as:
%The semi--ring structure and the spectral pairs of
%sesqui--linear forms, {\em Algebra Colloquium}, {\bf 1}:1, 1994, 85-95.
%
%\bibitem{[17]} N\'emethi, Andr\'as:
%Variation structures: results and open problems,
%{\em Banach Center Publications}, Volume {\bf 33}, 245-257;
% {\em Symposium on Singularities and
%Differential Equations}, Autumn 1993.
%
%\bibitem{RealS} N\'emethi, Andr\'as: The real Seifert form and
%the spectral pairs of isolated hypersurface singularities,
%{\em Compos. Math.}, {\bf 98}, 23-41, 1995.

%\bibitem{[19]} N\'emethi, A.:
%The equivariant signature of hypersurface singularities and
%eta--invariant, {\em  Topology} {\bf Vol. 34}, No. 2, 243-259, 1995.

%\bibitem{CR} N\'emethi, Andr\'as:
%The mixed Hodge structure of a complete intersection with
%isolated singularity, {\em C. R. Acad. Sci. Paris}, t.{\bf 321}, S\'erie I,
%447-452, 1995.

%\bibitem{[21]} N\'emethi, Andr\'as:
%The eta--invariant of variation structures. I.,
%{\em Topology and its Applications}, {\bf 67}, 95-111, 1995

%\bibitem{NTS} N\'emethi, Andr\'as:
%On the Seifert form at infinity associated with polynomial maps,
%{\em Journal of the Math. Soc. of Japan}, {\bf 51}, No 1, 63-70, 1999. 
%
%\bibitem{[28]} N\'emethi, Andr\'as:
%On the spectrum of curve singularities, 
%{\em Proceedings of the Singularity Conference}, Oberwolfach, July 1996;
%Progress in Mathematics, {\bf Vol. 162}, 93-102, Birkh\"auser 1998.
%
%\bibitem{[26]} N\'emethi, Andr\'as:
%Hypersurface singularities with 2--dimensional critical locus,
%will appear in {\em Journal of the London Math. Soc.}.
%
%\bibitem{[27]} N\'emethi, Andr\'as:
%Generalized Weil's Reciprocity Law and Multiplicativity Theorems,
%{\em Transaction of the AMS}, {\bf 349}, 2687-2697,  1997.

%\bibitem{[29]} N\'emethi, A.: The signature of $f(x,y)+z^n$, 
%{\em Proceedings of Real and Complex Singularities}, Liverpool,  August 1996; 
%London Math. Soc. Lecture Notes Series, {\bf 263}, 131-149, 1999.

%\bibitem{[30]} N\'emethi, Andr\'as:
%Some topological invariants of isolated hypersurface
%singularities, Five lectures of the EMS--Summer School, Eger (Hungary)
%1996, will appear in the {\em Proceedings of the Summer school}.
%\ http://www.mi.aau.dk/\~\,esn/titles\_1997.html

\bibitem{Five} N\'emethi, A.: Five lectures on normal surface singularities,
lectures delivered at the Summer School in {\em Low dimensional topology}
Budapest, Hungary, 1998; Bolyai Society Math. Studies {\bf 8} (1999), 
269-351.

%\bibitem{[35]} N\'emethi, A.: 
%Casson invariant of cyclic coverings via eta--invariant and
%Dedekind sums,  {\em Topology and its Applications}, {\bf 102}(2), 181-193, 
%2000.

\bibitem{[36]} N\'emethi, A.:
Dedekind sums and the signature of $f(x,y)+z^N$, 
{\em Selecta Mathematica}, New series, {\bf 4}, 361-376, 1998. 

\bibitem{[37]} N\'emethi, A.:
Dedekind sums and the signature of $f(x,y)+z^N$,II.,  
{\em Selecta Mathematica}, New series, {\bf 5}, 161-179, 1999. 

\bibitem{Ninv}  N\'emethi, A.:
``Weakly'' Elliptic Gorenstein Singularities of
Surfaces, {\em Inventiones math.}, {\bf 137}, 145-167, 1999. 

%\bibitem{[34]} N\'emethi, A.:
%On the Birkhoff normal form of a completely integrable
%Hamiltonian system near a fixed point with resonance (joint paper
%with T. Kappeler and Y. Kodama), 
%{\em Annali della Scuola superiore de Pisa}, {\bf 26} (4), 623-661, 1998.

%\bibitem{NSz1} N\'emethi, A.:
%The resolution of some surface singularities, I., 
%(cyclic coverings);  
% Proceedings of the AMS
%Conference, San Antonio, 1999; 
%{\em Contemporary Mathematics} {\bf 266}, 89-128.
%Singularities in Algebraic 
%and Analytic Geometry, (C. G. Melles and R. I. Michler Editors), American
%Math. Soc. 2000, 

\bibitem{NOSZ} N\'emethi, A.: On the Ozsv\'ath-Szab\'o invariant of negative 
definite plumbed 3-manifolds, {\sf arXiv:math.AG/0310083}.

\bibitem{Line} N\'emethi, A.: Line bundles associated with normal surface 
singularities, {\sf arXiv:math.AG/0310084}.

\bibitem{INV} N\'emethi, A.: Invariants of normal surface singularities,
to appear in 
the Proceedings of the Conference: {\em Real and Complex Singularities}, 
San Carlos, Brazil, August 2002. 

\bibitem{[51]} N\'emethi, A. and Nicolaescu, L.I.: 
 Seiberg-Witten invariants and surface singularities,
{\em Geometry and Topology},  Volume {\bf 6} (2002), 269-328. 

\bibitem{[52]} N\'emethi, A. and Nicolaescu, L.I.: 
 Seiberg-Witten invariants and surface singularities II
(singularities with good ${\bf C}^*$-action), {\sf arXiv:math.AG/0201120;}
to appear in the {\em Journal of LMS}. 

\bibitem{[55]} N\'emethi, A.  and Nicolaescu, L.I.: 
 Seiberg-Witten invariants and surface singularities III
(splicings and cyclic covers), {\sf  arXiv:math.AG/0207018.}

%\bibitem{NS}  N\'emethi, Andr\'as and Sabbah, Claude:
%Semicontinuity of the spectrum at infinity, 
%to appear in {\em Abh. Math. Seminar der Universit\"at Hamburg}. 
%http://math.polytechnique.fr/cmat/sabbah/articles.html

%\bibitem{[25]} N\'emethi, Andr\'as and Steenbrink, Joseph:
%On the monodromy of curve singularities,
%{\em Math. Zeitschrift}, {\bf 223}, 587-593, 1996.
%
%\bibitem{Ann}  N\'emethi, Andr\'as and Steenbrink, Joseph:
%Extending Hodge bundles for Abelian Variations,
%{\em Annals of Math.}, {\bf 143}, 131-148, 1995.
%
%\bibitem{NYJ} N\'emethi, Andr\'as and Steenbrink, Joseph:
%Spectral pairs, mixed Hodge modules and series of plane
%curve singularities, {\em New York Journal of Mathematics}, August 16, 1995;
%(http://nyjm.albany.edu:8000/j/v1/Nemethi-Steenbrink.html)

%\bibitem{NSz2} N\'emethi, Andr\'as and Szil\'ard, \'Agnes: 
%The resolution of some surface singularities, II.,
%(Iomdin's series),  Proceedings of the AMS
%Conference, San Antonio, 1999;
%{\em Contemporary Mathematics} {\bf 266}, Singularities in Algebraic 
%and Analytic Geometry, (C. G. Melles and R. I. Michler Editors), American
%Math. Soc. 2000, 129-164.

\bibitem{Neu} Neumann, W.:  Abelian covers of quasihomogeneous surface
singularities, {\em Proc. of Symposia in
Pure Mathematics}, vol. 40, Part 2, 233-244.

%\bibitem{Np} Neumann, Walter D.: Cyclic suspensions of knots and periodicity
%of signature for singularities, {\em Bull. Amer. Math. Soc.},
% Volume {\bf 80}, Number 5, 977-981, 1974.

%\bibitem{Ntr} Neumann, Walter D.: 
%Signature related invariants of manifolds -- I.
%Monodromy and $\gamma$--invariants, {\em Topology}, Vol. {\bf 18}, 147-172,
%1979.

%\bibitem{Neu1} Neumann, Walter D.: Invariants of plane curve singularities,
%{\em Monographie No 31 de L'Enseignement Math\'ematique}, 223-232, 1983.

%\bibitem{Neu2} Neumann, Walter D.: 
%Splicing Algebraic Links, {\em Advenced Studies
%in Pure Math.}, {\bf 8}, Complex Analytic singularities, 349-361,  1986.

%\bibitem{NP} Neumann, W.D.: A calculus for plumbing applied to the
%topology of complex surface singularities and degenerating complex curves.
%{\em Transactions of the AMS}, {\bf 268} Number 2, 299-344, 1981.

%\bibitem{NeR} Neumann, W.D.  and Raymond, F.: Seifert manifolds, plumbing,
% $\mu$-invariant and orientation reserving maps, {\em Algebraic and Geometric
%Topology},
% (Proceedings, Santa Barbara 1977), 
%Lecture Notes in Math. {\bf 664}, 161-196.

\bibitem{NW} Neumann, W. and Wahl, J.: Casson invariant of links
of singularities, {\em Comment. Math. Helv.} {\bf 65}, 58-78, 1991.

\bibitem{NWnew} Neumann, W. and Wahl, J.:
Universal abelian covers of surface singularities,
{\em Trends in singularities}, 181-191, Trends Math.,
Birkh\"auser, Basel, 2002.

\bibitem{NWnew2} Neumann, W. and Wahl, J.:
Universal abelian covers of quotient-cusps, {\em Math. Ann.} {\bf 326}, 75-93, 2003.

\bibitem{NWuj} Neumann, W. and Wahl, J.: 
Complex surface singularities with integral homology sphere links,
{\sf arXiv:math.AG/0301165}.

%\bibitem{Nico5} Nicolaescu, L.I.: 
%Seiberg-Witten invariants of rational homology  spheres, 
%{\sf  math.DG/0103020.}

\bibitem{Ok} Okuma, T.: Numerical Gorenstein elliptic singularities,
preprint 

%\bibitem{O} Ono, Isao and Watanabe, Kimio: On the singularity of 
%$z^p+y^q+x^{pq}=0$,
%{\em Sci. Rep. Tokyo Kyoika Daigaku Sect. A}, {\bf 12}, 123-128, 1974.

%\bibitem{OW} Orlik, P. and Wagreich, Ph.: Isolated singularities of
%algebraic surfaces with $\bc^*$ action, {\em Ann. of Math.} (2) {\bf 93},
%205-228, 1971.

%\bibitem{Osh}  O'Shea, D.:  Topological trivial deformations of isolated
%quasi-homogeneous hypersurface are equimultiple, {\em Proc. AMS}, 
%{\bf 101} (1987) 260-262.

%\bibitem{OSz} Ozsv\'ath, P.S. and Szab\'o, Z.: Holomorphic discs
%and topological invariants for rational homology three-spheres,
%math.SG/0101206.

%\bibitem{OSzP} Ozsv\'ath, P.S. and Szab\'o, Z.: On the Floer 
%homology of plumbed three-manifolds, math.SG/0203265.

%\bibitem{Peron} Perron, B.: $\mu$ ``constant'' implique ``type topologic 
%constant'' en dimension troi, Universite de Dijon, preprint ???????????

%\bibitem{Pic1} Pichon, A.:  Singularities of Complex  surfaces  with
%Equations $z^k=f(x,y)=0$, {\em International Math. Research Notices}, {\bf 5} 
%(1997), 241-246.

%\bibitem{Pi} Pinkham, H.: Deformations of algebraic varieties 
%with $G_m$ action, {\em Ast\'erisque} {\bf 20} (1974), 1-131.

\bibitem{Pi1}  Pinkham, H.: Normal surface singularities with 
${\bf C}^*$--action, {\em Math. Ann.,} {\bf 227}, 183-193,  1977.

%\bibitem{Pi2} Pinkham, Henry: Deformations of Normal Surface Singularities
%with ${\bf C}^*$ Action, {\em  Math. Ann.}, {\bf 232}, 65-84,  1978.

%\bibitem{P1}Pinkham, H.: Smoothing of the $D_{pqr}$ singularities,
%$p+q+r=22$,  {\em Proc. of Symp. in Pure Math.}, {\bf 40 } Part 2, (1983),
%373-377.

%\bibitem{Po} Pommersheim, James E.: Toric varieties, lattice points and
%Dedekind sums, {\em Math. Ann}, {\bf 295}, 1-24,  1993.


%\bibitem{Ra} {\bf H Rademacher}, {\it Some remarks on certain generalized
%Dedekind sums}, Acta Arithmetica, {9}(1964), 97-105.

%\bibitem{RG} Rademacher, H. and Grosswald, E.: Dedekind Sums, 
%The Carus Math. Monographs, MAA, 1972.

%\bibitem{Rad} Rademacher, Hans: Generalization of the Reciprocity formula for
%Dedekind sums, {\em Duke Math. Journal}, {\bf 21}, 391-397,  1954.

%\bibitem{Ratiu} Ratiu, A.: The Jones--Witten invariants of tree manifolds,
%These de Doctorat, Universite Paris 7, 1996.

%\bibitem{MR}  Reid, M.: Chapters on Algebraic Surfaces.
%In: Complex Algebraic Geometry,
%IAS/Park City Mathematical Series,  Volume {\bf 3}  (J. Koll\'ar editor),
%3-159, 1997.

%\bibitem{MRCo} Reid, Miles: Canonical 3-folds, Algebraic Geomtry Angers,
%1979 (A. Beauville, editor) 273-310; Sijthoff \& Noordhoff, 1980.

%\bibitem{Riem} Riemenschneider, O.: Bemerkungen zur Deformationstheorie 
%Nichtrationaler Singularit\"aten, {\em Manus. Math.}, {\bf 14}, 91-99, 1974.


%\bibitem{Ri} Riemenschneider, O.: Deformationen von Quotintensingularit\"aten
%(nach zyklischen Gruppen), {\em Math. Ann.}, {\bf 209} (1974), 211-248. 

%\bibitem{saeki} Saeki, O.: Topological types of complex isolated 
%hypersurface singularities, {\em Kodai Math. J.}, {\bf 12} (1990), 23-29.

%\bibitem{saeki2} Saeki, O.: Real Seifert form  determines the spectrum for
%quasihomogeneous hypersurface singularities in $\bc^3$, {\em J. Math. Soc. 
%Japan}, {\bf 52} (2) (2000), 409-431. 

%\bibitem{KSa} Saito, Kyoji: The Zeroes of Characteristic Function $\chi_f$
%for the Exponents of a Hypersurface Isolated Singular Point,
%{\em Advanced Studies in Pure Math.} {\bf 1}, Algebraic Varieties and Analytic
%Varieties, 195-217, 1983.

%\bibitem{SaP} Saito, Morihiko: On the exponents and the geometric genus of an
%isolated hypersurface singularity, {\em Proc. of Sympos. in Pure Math.},
%{\bf 40}, Part 2, 465-472, 1983.

%\bibitem{Sa} Sakamoto, Koichi: 
%The Seifert matrices of Milnor fiberings defined by
%holomorphic functions, {\em J. Math. Soc. Japan}, {\bf 26}, 714-721, 1974.

%\bibitem{Scha} Schaps, M.: Deformations 
%of Cohen-Macauley Schemes of codimension
%2 and Non-Singular Deformations of Space Curves, {\em Am. J. Math.}, {\bf 
%99} (1977), 669-685.  


%\bibitem{SS} Scherk, John and Steenbrink, Joseph H. M.: On the Mixed Hodge
%Structure on the Cohomology of the Milnor Fiber, {\em Math. Ann.},
%{\bf 271}, 641-665, 1985.

%\bibitem{Sch} Schlessinger, M.: Functors of Artin Rings, {\em Trans. AMS},
%{\bf 130} (1968), 208-222. 

%\bibitem{SSS} Schrauwen, R. and Steenbrink, J. and Stevens, J.: Spectral
%Pairs and Topology of Curve Singularities, {\em Proc. Sumpos. Pure Math.},
%{\bf 53}, 305-328, 1991.

%\bibitem{Seade} Seade, J.A.:  A cobordism invariant for surface
%singularities, {\em Proc. of Symp. in Pure Math.}, {\bf 40} (2) (1983),
%479-484.

%\bibitem{ST} Sebastiani, Marcos and Thom, Ren\'e: 
%Un r\'esultat sur la monodromie,
%{\em Inventiones Math.}, {\bf 13}, 90-96, 1971.

\bibitem{GPS01} G.-M. Greuel, G. Pfister, and H. Sch\"onemann. 
{\sc Singular} 2.0. A Computer Algebra System for Polynomial
Computations. Centre for Computer Algebra, University of
Kaiserslautern (2001). {\tt http://www.singular.uni-kl.de}.

%\bibitem{Spiv} Spivakovsky, M.: Sandwiched singularities and desingularization
%of surfaces by normalized Nash transformations, {\em Annals of Math.}, 
%{\bf 131} (1990), 411-491.

%\bibitem{StO} Steenbrink, J.H.M.: Mixed Hodge structures on the 
%vanishing
%cohomology, Nordic Summer School/NAVF, Symposium in Mathematics, Oslo, 1976.

%\bibitem{St.Proc} Steenbrink, J.H.M.: Mixed Hodge structures
%associated with isolated singularities, {\em Proc.  Symp.
%Pure Math.}, {\bf 40}, Part 2  (1983), 513-536.

%\bibitem{St1} Stevens, J.: Elliptic Surface
%Singularities and Smoothings of Curves,
%{\em Math. Ann.}, {\bf 267}, 239-247,  1984.

%\bibitem{St2} Stevens, J.: Kulikov singularities; Thesis, Leiden 1985

%\bibitem{stev} Stevens, J.: On the versal deformation of cyclic 
%quotient singularities, {\em LNM}, {\bf 1462} (1991), 302-319.
%(Singularity theory and its applications, Warwick 1989)

%\bibitem{stev2} Stevens, J.: Partial resolutions of rational quadruple
%points, {\em Int. J. of Math.}, {\bf 2} (2) (1991), 205-221.

%\bibitem{Cyc} Teissier, B.: Cycles evanescents, sections planes et
%conditions de Whitney, {\em Asterisque} {\bf 7-8} (1973), 285-362.

%\bibitem{teissier} Teissier, B.: Deformation \`a type topologique constant
%II, {\em S\'eminaire Douady-Verdier 1972}. 

%\bibitem{IaII} Teissier, B.: R\'esolution simultan\'ee I, II, 
%{\em LNM} {\bf 777} (1980), 71-146.

%\bibitem{Tju} Tjurina, G.-N.: Locally Flat Deformations of Isolated 
%Singularities of Complex Spaces, {\em Math. USSR Izvestia} {\bf 3} (1969), 
%967-999.

%\bibitem{To} Tomari, Masataka: The inequality $8p_g\leq \mu$ for hypersurface
%two--dimensional isolated double points, {\em Math. Nachr.}, {\bf 164},
%37-48, 1993.

\bibitem{To1} Tomari, M.: A $p_g$--formula and elliptic singularities,
Publ. R.I.M.S. Kyoto University, {\bf 21} (1985), 297-354.

%\bibitem{To2} Tomari, Masataka and Watanabe, Kei-ichi: 
%Filtered rings, Filtered Blowing--Ups,
%Normal Two--Dimensional Singularities with ``Star--Shaped'' Resolution,
%Publ. R. I. M. S. Kyoto University,
%{\bf 25}, 681-740, 1989.

\bibitem{Tu5} Turaev, V.G.:   Torsion invariants of $Spin^c$-structures 
on $3$-manifolds, {\em Math. Res. Letters}, {\bf 4} (1997), 679-695.

%\bibitem{Va} Vaqui\'e, M.: R\'esolution  simultan\'ee de surfaces normales,
%{\em Ann. Inst. Fourier}, {\bf 35} (1985), 1-38.

%\bibitem{V} Varchenko, Alexander: Asymptotic Hodge structure in the vanishing
%cohomology, {\em Math. USSR Izv.}, {\bf 18}, 469-512, 1982.

%\bibitem{Wagreich} Wagreich, Ph.: Elliptic singularities of surfaces.
%Amer. J. of Math., {\bf 92}, 419-454, 1970.


%\bibitem{Wag} Wagreich, Ph.: {\sl The structure of quasihomogeneous
%singularities}, Proc. Symp. in Pure Math., {\bf 40}:2, 598-611 (1983)

%\bibitem{Walker} Walker, K.: An extension of the Casson's  invariant,
%Ann. of Math. Studies {\bf 126}, Princeton University Press, 1992.

%\bibitem{Whal} Wahl, M.J.: Vanishing Theorems for Resolutions of
%Surface Singularities. Inventiones math., {\bf 31}, 17-41, 1975.

\bibitem{WahlAnnals} Wahl, M.J.: Equisingular deformations of normal 
surface  singularities, I,  {\em Ann. of Math.}, {\bf 104}, 325-356, 1976.

%\bibitem{WComp} Wahl, M. J.: Simultaneous resolution of rational 
%singularities, {\em Compositio Math.}, {\bf 38} (1) (1979), 43-54.

%\bibitem{WEll} Wahl, M. J.: Elliptic Deformations of Minimally Elliptic
%Singularities, {\em Math. Ann.}, {\bf 253} (1980), 241-262.

%\bibitem{W} Wahl, J.: Smoothings of normal surface singularities,
%{\em Topology}, {\bf 20}, 219-246, 1981.

%\bibitem{Wat} Watanabe, K.: On Plurigenera of Normal Isolated 
%Singularities I, {\em Math. Ann.}, {\bf 25} (1980), 65-94. 

%\bibitem{Weber} Weber, C.: Problems, {\em Noeuds, tresses  et singularit\'es},
%Monographie {\bf 31} de L'Enseignement Math. (1983), 247-259.

%\bibitem{XY1} Xu, Y. and Yau, S.S.-T.: The inequality $\mu \geq 12p_g-4$
%for hypersurface weakly elliptic singularities, {\em Contemporary
%Math.}, {\bf 90}, 317-344, 1989; (Richard Randell Editor).

%\bibitem{XY2} Xu, Y. and Yau, Stephen S.-T.: Durfee's conjecture and 
%coordinate free
%characterization of homogeneous singularities, {\em J. Diff. Geometry
%} {\bf 37}, 375-396, 1993.

%\bibitem{XY3} Xu, Y. and Yau, Stephen S.-T.: A sharp estimate of the number 
%of integral points in a tetrahedron, {\em J. reine angew. Math.}, {\bf 423},
%199-219, 1992.

%\bibitem{XY}  Xu, Y and Yau, S. S.-T. Yau:
% Classification of topological types 
%of isolated quasi-homogeneous two dimensional hypersurface singularities,
%{\em Manuscripta Math.}, {\bf 64} (1989), 445-469.

%\bibitem{YauConj} Yau, S.S.-T.: 
%Topological types of isolated hypersurface 
%singularities, {\em Contemporary Math.}, {\bf 101} (1989), 303-321.

\bibitem{Yau4} Yau, S.S.-T.: On almost minimally elliptic singularities,
{\em Bulletin of the AMS}, {\bf 83} Number 3, 362-364, 1977.

%\bibitem{Yau6} Yau, Stephen S.-T.: Normal singularities of surfaces,
%{\em Proc. of Symposia  in Pure Math.}, {\bf 32}, 195-198, 1978.

%\bibitem{Yau7} Yau, Stephen S.-T.:
% Normal two--dimensional elliptic singularities,
%{\em Transactions of AMS}, {\bf 254}, 117-134, 1979.

%\bibitem{Yau2} Yau, Stephen S.-T.:
% Hypersurface dual graphs of normal singularities
%of surfaces, {\em  Amer. J. of Math.}, {\bf 101}, 781-812, 1979.

%\bibitem{Yau3} Yau, S.S.-T.:
% Gorenstein singularities with geometric genus
%equal to one, {\em  Amer. J. of Math.}, {\bf 101}, 813-854, 1979.

\bibitem{Yau5} Yau, S.S.-T.: On strongly elliptic singularities,
{\em Amer. J. of Math.}, {\bf 101}, 855-884, 1979.

\bibitem{Yau1} Yau, S.S.-T.: On maximally elliptic singularities,
{\em Transactions of the AMS}, {\bf 257} Number 2, 269-329, 1980.

%\bibitem{Y1} Yoshihara, H.: Rational curve with one cusp,
%{\em Proc. AMS}, {\bf 89}(1) (1983), 24-26.

%\bibitem{Z} Zagier, Don: Higher dimensional Dedekind sums, {\em
%Math. Ann.}, {\bf 202}, 149-172, 1973.

%\bibitem{cangr} Zariski, Oskar: The reduction of singularities of an
%algebraic surface, {\em Ann. Math.}, {\bf 40} (2), 639-689, 1939.

\bibitem{Zariskiconj} Zariski, O.: Some open questions in the theory
of singularities, {\em Bulletin of the AMS}, {\bf 77}, 481-491, 1971.


\end{thebibliography}
\end{document}